\documentclass[11pt]{article}
\usepackage[utf8]{inputenc}
\usepackage[T1]{fontenc}
\usepackage[usenames,dvipsnames]{xcolor}
\usepackage[pdftex]{graphicx} 
\usepackage[english]{babel}  
\usepackage{fullpage} 
\usepackage{tikz} 
\usepackage{amsmath}
\usepackage{amssymb}
\usepackage{amsthm}
\usepackage{bbm}
\usepackage{cancel} 
\usepackage[shortlabels]{enumitem} 
\usepackage{amsfonts}
\usepackage{dsfont}

\usepackage{float}
\usepackage{subcaption}
\usepackage{comment}
\usepackage{standalone}

\newcommand{\prob}{\mathbb{P}}
\newcommand{\E}[1]{E\left[#1\right]}

\newcommand{\Exp}[1]{\mathbb{E}\left[#1\right]}

\newcommand{\Prob}[1]{\mathbb{P}\left(#1\right)}
\newcommand{\Var}{\operatorname{Var}}
\newcommand{\Cov}{\operatorname{Cov}}
\newcommand{\Expn}[1]{\mathbb{E}_0\left[#1\right]}
\newcommand{\Expa}[1]{\mathbb{E}_1\left[#1\right]}
\newcommand{\Varn}{\operatorname{Var}_0}
\newcommand{\Vara}{\operatorname{Var}_1}

\newcommand{\ind}[1]{\mathbbm{1}_{\{#1\}}}
\DeclareOldFontCommand{\rm}{\normalfont\rmfamily}{\mathrm}
\newcommand{\dd}{{\rm d}}

\newcommand{\bx}{\operatorname{x}}

\newcommand{\Wi}{W}
\newcommand{\plim}{\ensuremath{\stackrel{\prob}{\longrightarrow}}}

\newcommand{\clara}[1]{\color{blue}Clara: {#1} \color{black}}
\newcommand{\mitch}[1]{\color{ForestGreen}{#1}\color{black}}

\newcommand{\ER}{Erd\H{o}s-R\'enyi }

\numberwithin{equation}{section}

\begin{document}

\theoremstyle{plain}
    \newtheorem{thm}{Theorem}[section]
    \newtheorem{lemma}[thm]{Lemma}
    \newtheorem{cor}[thm]{Corollary}
	\newtheorem{prop}[thm]{Proposition}
	\newtheorem{assumption}{Assumption}
	
\theoremstyle{definition}
    \newtheorem{defn}{Definition}
    \newtheorem*{propr}{Property}
	
\theoremstyle{remark}
    \newtheorem*{remark}{Remark} 
    \newtheorem*{ex}{Example}

\title{Localized geometry detection in scale-free random graphs}
\author{Gianmarco Bet, Riccardo Michielan, Clara Stegehuis}

\maketitle

\begin{abstract}
We consider the problem of detecting whether a power-law inhomogeneous random graph contains a geometric community, and we frame this as an hypothesis testing problem. More precisely, we assume that we are given a sample from an unknown distribution on the space of graphs on $n$ vertices. Under the null hypothesis, the sample originates from the inhomogeneous random graph with a heavy-tailed degree sequence. Under the alternative hypothesis, $k=o(n)$ vertices are given spatial locations and connect between each other following the \textit{geometric} inhomogeneous random graph connection rule. The remaining $n-k$ vertices follow the inhomogeneous random graph connection rule. We propose a simple and efficient test, which is based on counting normalized triangles, to differentiate between the two hypotheses. We prove that our test correctly detects the presence of the community with high probability as $n\to\infty$, and identifies large-degree vertices of the community with high probability.
\end{abstract}

\section{Introduction}
Random graphs provide a unified framework to model many complex systems in biology, computer science, sociology, as well as numerous other sciences. The random graph paradigm usually involves specifying a probabilistic mechanism to generate a graph, and then studying the properties of the resulting network, be it topological (connectedness, clustering, etc) or statistical (fit to data). Random graphs are particularly useful as \textit{null models} to determine if some observed real-world network deviates from its expected structure in a statistically significant way. In this context, it has been widely observed that real-world networks share two defining features: heavy-tailed degree sequences and large clustering~\cite{faloutsos1999,newman2003b}. Both these features are not reproduced by the classical \ER random graph model, which makes this an unsatisfactory null model for most applications. Consequently, alternative models have been developed to match the degree sequence and clustering observed in real-world networks. The so-called \textit{inhomogeneous random graph} (IRG)~\cite{chung2002} is a popular generalization of the \ER random graph obtained by assigning weights to nodes, and connecting two nodes with a probability that is proportional to the products of their weights. This way, the IRG can reproduce an arbitrary degree sequence, but still has low clustering. 

A popular method to obtain a model with a large clustering is to embed the vertices in a metric space (such as the sphere or the torus) and connecting them with probabilities proportional to their distances~\cite{krioukov2010}. Indeed, the presence of distances makes two neighbors of a given vertex likely to be close by and therefore connected as well, due to the triangle inequality.  By embedding the vertices of the IRG in a torus, one obtains the so-called \textit{geometric inhomogeneous random graph} (GIRG)~\cite{bringmann2019}. This model creates networks with two phenomena that are often observed in real-world networks: heavy-tailed degree sequences, as well as high clustering. However, it is often the case that clustered nodes are not spread evenly across the network, but rather they form communities. It is then of great practical interest, first, to establish if these communities are actually present, and, second, to identify them. Perhaps surprisingly, the early literature on the latter problem did not address the former \cite{girvan2002community,newman2006modularity,newman2004finding}. In fact, often the focus of the community detection literature lies on algorithms to extract a community structure from given networks, regardless of whether the structure is actually present. These algorithms are usually tested on random graph models with a known community structure. One such example is the stochastic block model~\cite{holland1983}, which has received considerable attention due to its mathematical tractability. However, this comes at the expense of unrealistic assumptions, such as very large communities (of the same order of the graph size) and homogeneous degree distributions. To overcome this, Arias-Castro and Verzelen \cite{arias2014community,verzelen2015community} considered the problem of detecting the presence of a small community in an \ER random graph. They find the region in the parameter space where (almost sure) detection is impossible, and they give tests that are able to detect the community outside of this region. These results were later generalized to the IRG in \cite{bogerd2021detecting}. However, the results in \cite{bogerd2021detecting} do not apply to heavy-tailed degree sequences and the community is obtained by tuning an ad-hoc density parameter. In this paper, we take a further step towards obtaining detection results for realistic networks. More precisely, our null model is the IRG with an heavy-tailed degree sequence, and the community, if present, is obtained by embedding a small number of the total nodes within a torus, and connecting them according to the usual GIRG connection probabilities. Due to the geometric nature of the community, it contains many triangles, as opposed to the tree-like nature of an IRG-based community. This realistic feature of the community allows us to develop efficient testing methods for the testing as well as the identification of the community.

More precisely, when the community is indeed present in a graph of size $n$, $n-k$ nodes of the network form connections with each other on a non-geometric basis. The other $k$ nodes have a position in some geometric space, and nearby nodes are more likely to connect. This geometric setting creates a subgraph with many triangles, and can therefore be thought of as a community in the network. This geometric structure is less restrictive than planting a clique, and more realistic than a dense inhomogeneous random graph as a community. Furthermore, the fact that the planted structure is geometric allows for efficient, triangle-based tests to detect and identify the structure. To the best of our knowledge, this type of planted structure has not been considered before. 

Our contributions are the following:
\begin{itemize}
    \item We provide a statistical test to detect the presence of a planted geometric community. Unlike other detection tests for the purpose of dense subgraph detection~\cite{bogerd2021detecting,verzelen2015community}, the test works for heterogeneous degrees. This method is triangle-based, making it efficient in implementation. Rather than using standard triangle counts, which may not be able to differentiate between geometric and non-geometric networks, the statistic weighs the triangles based on their evidence for geometry. 
    \item We provide a statistical method to identify the largest-degree vertices of the planted geometric community. This method is also triangle-based, and therefore it is efficient to implement. We show that this method achieves exact recovery among all high-degree vertices.
    \item We provide a method to infer the size of the planted geometric community
    This method uses the largest-degree identified vertices of the planted geometric community to obtain an estimate for the community size 
    based on the convergence of order statistics. 
    \item We show numerically that the combination of these tests leads to an accurate identification of the planted geometric community. Furthermore, these tests can be performed in computationally only $O(n^{3/2})$ time~\cite{latapy2007practical}.
\end{itemize}



\subsection{Related literature}
Our work lies at the intersection of two rich lines of research: community detection and what might be referred to as \textit{structure detection}. In the former setting, one is given a sample from a known random graph model and the task is to determine if there is a statistically-unlikely dense subgraph, and possibly to identify it. In this context, the \textit{planted clique problem} has received considerable attention as a testbed for community detection algorithms. In this model, a large network of size $n$ is generated according to some mechanism, and a small clique of size $k$ might be planted in it~\cite{alon1998finding,wu2021statistical}. The seminal works \cite{arias2014community,verzelen2015community} form a stepping stone towards more realistic dense communities. Their null model is the (resp.~dense, sparse) \ER random graph, and, when present, the community is a small subset of vertices with larger connection probability than in the null model. See also \cite{hajek2015computational}. Further generalizing this work, \cite{bogerd2021detecting} focuses on detecting a dense subgraph in an inhomogeneous random graph. More precisely, their null model is the inhomogeneous random graph, and in the alternative hypothesis the connection probabilities of a small subset of nodes $C$ are increased by a multiplicative factor $\rho_C>1$. Crucially, their approach requires a precise control on the inhomogeneity of the graph and does not work, for example, for heavy-tailed degree distributions. Therefore, the difference between our work and \cite{bogerd2021detecting} is two-fold. First, we consider the case of power-law vertex weights, which is more attractive from a modelling point of view. Second, the planted structure is a community by virtue of the underlying geometrical structure, rather than by tuning an additional model parameter of a tree-like graph. The work \cite{bet2021detecting} tackles the opposite problem to ours, namely detecting mean-field effects in a geometric random graph model. More precisely, their null model is a geometric random graph, and in the alternative hypothesis a small subset of vertices connects with every other vertex according to i.i.d.~Bernoulli random variables. They provide detection thresholds, as well as asymptotically powerful tests.

On the other hand, in the setting of structure detection, one is given a sample from an unknown random graph model, and the task is to determine if the sample originates from a mean-field model or a structured (e.g., geometric) model. In \cite{bubeck2016testing} (see also \cite{eldan2020information}), the null model is the \ER graph, and the alternative model is a high-dimensional geometric random graph. For recent progress on this problem, see \cite{brennan2020phase}. \cite{gao2017testing} proposes a test based on small subgraph counts to distinguish between the \ER graph and a general class of structured models which includes the stochastic block model and the configuration model. More recently, \cite{bresler2018optimal} proposes a test to distinguish between a mean-field model and Gibbs models, and~\cite{michielan2022} proposes a test to distinguish between a power-law random graph with and without geometry. See also \cite{ghoshdastidar2020two} for two-sample hypothesis testing for inhomogeneous random graphs. \cite{jin2015fast} proposes the so-called SCORE algorithm for community detection on the Degree-Corrected Block Model. One of their main ideas is overcome the statistical issues caused by the heterogeneity of the degree distribution by constructing test statistics that are properly normalized so as to cancel out the effects of vertex weights. In a similar spirit, \cite{jin2018network} considers an inhomogeneous random graph with community and proposes a normalized test based on short paths and short cycles to detect the presence of more than one community.  Our work here is also graphlet-based (triangles in this case), but rather than taking all triangles as equal, we weigh the triangles based on the inhomogeneity of the network degrees. This provides a robust statistic to infer communities in heavy-tailed networks. 





\subsection{Structure of the paper}
The rest of the paper is structured as follows. In Section~\ref{sec:model} we explain the model and the hypotheses for our tests. In Section~\ref{sec:results} we provide the tests that we propose, and state our main results on their accuracy, followed by a discussion in Section~\ref{sec:discussion}. We finally prove our main results on detecting the presence of a geometric structure in Section~\ref{sec:detect}, and our results on the identification of the geometric structure in Section~\ref{sec:identify}.

\paragraph{Notation.}
We adopt the standard notation of a statistical testing problem. The null hypothesis will be denoted by $H_0$, and the alternative hypothesis by $H_1$. When operating under $H_0$, that is, assuming the null hypothesis holds, the probability of some event $E$ will be denoted by $P_0(E) := \mathbb{P}(E \;|\; H_0)$. We denote the expected value and the variance with respect to this probability measure by $E_0$ and $\Var_0$ respectively. On the other hand, when $H_1$ is assumed to hold, we will similarly use the notation $P_1$, $E_1$, $\Var_1$.
Throughout the entire paper, we will make use of the standard Bachmann-Landau notation. We write $f(n) = o(g(n))$ if $\lim_{n \to \infty} f(n)/g(n) = 0$, $f(n) = O(g(n))$ if $\limsup_{n \to \infty} f(n)/g(n) < \infty$ and $f(n) = \Omega(g(n))$ if $g(n) = O(f(n))$.
Finally, we say that a sequence of events $\{E_n\}_{n \geq 1}$ happens with high probability (w.h.p.) if $\lim_{n \to \infty} \mathbb{P}(E_n) = 1$.

\section{Model}\label{sec:model}
We now formulate the problem of community detection in a graph as an hypothesis testing problem. We are given a single sample of a simple graph $G=(V,E)$, where $V=[n]:=\{1,\ldots,n\}$ is the set of nodes, and $E\subseteq\{(i,j)\in V\times V: i< j\}$ are the edges. Note that by assumption $G$ does not contain self-loops and multiple edges. 

\paragraph{Null model.}
Under the null hypothesis $H_0$, the graph $G$ is a sample of the inhomogeneous random graph (IRG) model, which is defined as follows~\cite{chung2002}. 
To each vertex $i\in V$ we assign a weight $w_i$, 
where
\begin{equation*}
    F_n(x) = \frac{1}{n}\sum_{i \in V} \mathbbm{1}_{\{w_i \leq x\}}
\end{equation*}
denotes the empirical weight distribution. $F_n$ can also be seen as the weight distribution of a uniformly chosen vertex in the graph. We require the weight sequence to satisfy the following assumption.

\begin{assumption}\label{ass:degrees}There exist $\tau \in (2,3)$ and $C, w_0 > 0$, such that for all $x \geq w_0$ with $x \leq n^{1/(\tau-1)}/ \log(n)$ 
\begin{equation*}
    1 - F_n(x) = C x^{1-\tau} (1 + o(1)).
\end{equation*}
\end{assumption}

\noindent Given the weight sequence $\{w_i\}_{i\in V}$, any edge $(i,j)$ is present with probability
\begin{equation}\label{eq:IRG_connection_rule}
    p_{ij} = p(w_i,w_j) := \min\left(\frac{w_i w_j}{\mu n}, 1\right),
\end{equation}
independently from all other edges, where $\mu=w_0 \frac{\tau-1}{\tau-2}$. In Lemma \ref{lemma:sumweights} we prove that when $n$ is large, $\mu$ is asymptotically equal to the average weight.

\paragraph{Alternative model.}
Under the alternative hypothesis $H_1$, $k$ of the vertices form a community. Without loss of generality, we assume these are $V_C := \{1,\ldots,k\}\subset V$. For convenience, we denote $V_I:= V\setminus V_C$, and we call the elements of $V_I$ \emph{type-A} vertices, while we call the elements of $V_C$ \emph{type-B} vertices. Let us now define the geometric community more precisely.
Let $\mathcal X=[0,1]^d$ be the $d$-dimensional torus. We endow $\mathcal X$ with the norm
\begin{align}
\|x-y\| = \sup_{i=1,\ldots,d}\min\{\vert x(i)-y(i)\vert, \vert1-(x(i)-y(i))\vert\}, 
\end{align}
where $x=(x(1),\ldots, x(d))$ and $y=(y(1),\ldots,y(d))$ are elements of $\mathcal X$. Note that this is the usual infinity norm compatible with the torus structure. To each vertex $i\in V_C$ we assign a (random) position $\bx_i$ in the torus $\mathcal{X}$. Formally, $(\bx_i)_{i\in V_C}$ is a sequence of random variables distributed uniformly over $\mathcal X$, and we will denote by $(x_i)_{i\in V_C}$ a realization of such random sequence. 
Again, we assign to each vertex $i\in V$ a weight $w_i$, where $(w_i)_{i \in V}$ is a sequence satisfying Assumption~\ref{ass:degrees}. Additionally, defining the empirical distribution of the vertex weights in the geometric community as
\begin{equation*}
    F_k(x) = \frac{1}{k}\sum_{i \in V_C} \mathbbm{1}_{\{w_i \leq x\}},
\end{equation*}
we will also require that $F_k(x)$ has a power law tail.
\begin{assumption}\label{ass:degreesgeo} Let $\tau, C, w_0$ be the same as in Assumption \ref{ass:degrees}. Then, for all $x \geq w_0$ with $x = O\left(k^{1/(\tau-1)}/ \log(k)\right)$,
\begin{equation*}
    1 - F_k(x) = C x^{1-\tau} (1 + o(1)).
\end{equation*}
\end{assumption}

\noindent Under $H_1$, any edge $(i,j)\in V_I\times V$ is present independently from all other edges with probability given by 
\begin{equation}\label{eq:corr_IRG_connection_rule}
    p_{ij} = p(w_i,w_j) := \frac{1}{1+C_1}\min\left(\frac{w_i w_j}{\mu n}, 1\right),
\end{equation}
where  $C_1 := (1 + (\gamma-1)^{-1})2^d$. That is, pairs with at least one type-A vertex connect with probability determined by the weights of the two endpoints, similarly as under $H_0$.
Moreover, any edge $(i,j)\in V_C\times V_C$ is present independently from all other edges with probability
\begin{equation}\label{eq:corr_GIRG_connection_rule}
    p_{ij} = p(w_i,w_j,x_i,x_j) := \frac{1}{1+C_1}\min\left(\frac{w_i w_j}{\mu k||x_i - x_j||^d}, 1\right)^\gamma.
\end{equation}
for some $\gamma \in (1,\infty]$. This is a 
geometric connection probability on $k$ vertices similarly to the GIRG model~\cite{bringmann2019}, multiplied by the factor $1/(1+C_1)\in(0,1)$. By convention, the choice $\gamma = \infty$ corresponds to the \textit{threshold} connection rule, that is, $p_{ij}=1/(1 + C_1)$ if $||x_i - x_j||^d \leq w_i w_j/(\mu k)$, and $p_{ij}=0$ otherwise. Thus, these $k$ type B vertices form connections based on their weights as well as their positions. In particular, the closer $x_i$ and $x_j$, the more likely they are to connect. The triangle inequality also ensures that a connection between type B nodes $i$ and $j$ and $i$ and $k$ makes it more likely for an edge between $j$ and $k$ to be present as well. Thus, the type B vertices are likely to be more clustered than the type A vertices. Note that an alternative interpretation for the connection rule \eqref{eq:corr_GIRG_connection_rule} is that it is the GIRG connection probability on $n$ vertices~\cite{bringmann2019}, where the positions of the vertices $V_C$ are sampled uniformly over the (shrinking) torus $[0,\lceil k/n\rceil]^d$. 

\paragraph{Sources of randomness.}
Observe that under $H_1$, two sources of randomness are present: the position sequence $(x_i)_{i \in V_C}$, and the random independent connections between vertices. Given a network sample, the positions that generated the network community are usually unknown when only observing the network connections. Thus, we assume that we have no knowledge of the positional vectors of the community. However, when a given network is realization of an inhomogeneous random graph, the degree of a vertex in the network is close to its weight, with high probability~\cite[Appendix C]{stegehuis2017}. In the case of the geometric inhomogeneous random graph, the same holds for all vertices with weights asymptotically larger than $\log^2(n)$ \cite{bringmann2019}. Therefore, in our setting it is reasonable to assume the weight sequence is known, as it would be possible to infer it from a the degree distribution of a given network.

\paragraph{Correction factor.} In the IRG, any vertex $i$ has expected degree $w_i(1+o(1))$. On the other hand, the random graph formed under $H_1$, without the correction factor $1/(1 + C_1)$ in \eqref{eq:corr_IRG_connection_rule} and \eqref{eq:corr_GIRG_connection_rule}, would introduce a bias on the expected power-law degree distribution. Therefore, a simple check on the degree distribution would be sufficient to determine if a random graph has been sampled from $H_0$ or $H_1$. With the correction factor of $1/(1+C_1)$, the expected degree of any vertex is $w_i (1 + o(1))$ under $H_1$ as well, as proved in Appendix \ref{app:consistency_correction}, which excludes a trivial detection test.

\section{Main results}\label{sec:results}
In this section we describe our main results regarding the detection and the identification of the geometric community. First, let us introduce a few important notions. A \textit{test} 
$\psi$ is a mapping from $G$ to $\{0, 1\}$. Here $\psi(G)=1$ indicates that the null hypothesis $H_0$ is rejected and the graph contains a planted geometric community, and $\psi(G) = 0$ otherwise. The \textit{risk} of such a test is defined as
\begin{equation}
R(\psi) := P_0(\psi(G)= 0) + P_1(\psi(G)  = 1).
\end{equation}%
%
Our goal is to distinguish $H_0$ and $H_1$ when the graph size $n$ is large. Formally, a sequence of tests $(\psi_n)_{n\geq 1}$ is said to be \textit{asymptotically powerful} when it has vanishing risk, that is $\lim_{n\to\infty }R(\psi_n)=0$. Such a sequence of tests identifies the underlying model correctly in the limit of $n\to\infty$.

\subsection{Detection}
In this section we first describe an asymptotically powerful test for planted geometric community detection, the \emph{weighted triangle test}. We will use the short-hand notation $\{i,j,k\}=\triangle$ to mean $\{(i,j), (j,k), (k,i)\}\subseteq E$. The test uses the \textit{weighted triangles} statistic
\begin{equation}
    W(G) := \sum_{a,b,c \in V} \frac{1}{w_a w_b w_c} \mathbbm{1}_{\{\{a,b,c\} = \triangle\}} \label{eq:W}.
\end{equation}
Thus, each triangle is given a weight that is inversely proportional to the product of the weights of its vertices. In this way, $W$ \textit{discounts} the triangles formed by high-weight vertices. Indeed, triangles between high-weight vertices are likely to be formed in geometric as well as in non-geometric random graphs. Therefore, standard triangle counts are not even able to distinguish between power-law geometric graphs and inhomogeneous random graphs~\cite{michielan2022}, and we need more advanced triangle-based statistics. The main distinction is given by the triangles formed between low-degree vertices, which are unlikely in non-geometric models. The weighted triangle test rejects $H_0$ when $W(G)$ is larger than some threshold $f(n)$. 
Formally, the weighted triangle test $\psi_W$ is defined as
\begin{align}\label{eq:testdetect}
    \psi_W(G) = \mathds 1_{\{W(G) \geq f(n)\}}.
\end{align}
The next result shows that there is significant freedom in the choice of $f(n)$, while still having an asymptotically powerful test: 

\begin{thm}\label{thm:detect}
    Let $f(n)$ be a function such that  $f(n)\to\infty$ as $n\to\infty$ and $f(n)= o(k)$. Then, the weighted triangle test is asymptotically powerful. 
\end{thm}

Theorem~\ref{thm:detect} shows that it is possible to detect the presence of any geometric subset as long as it grows (arbitrarily slowly) in $n$. Still, the test statistic ~\eqref{eq:testdetect} relies on knowledge of a lower bound on the geometric size $k$ when choosing the threshold $f(n)$, as Theorem~\ref{thm:detect} requires $f(n)=o(k)$.

%
\subsection{Identification}
We now focus on the problem of identifying the geometric vertices under $H_1$. 
When a test rejects $H_0$, the following goal is to identify the vertices that are part of the planted geometric part. To this end, let $\hat{V}_C\subseteq V$ be an estimator for the set of geometric vertices. We assume that the size of the planted geometric community, $k$ is known. 
To measure the performance of an estimator of the geometric vertices we use the risk function
\begin{equation}
    R_{id}(\hat{V}_C) := \mathbb{E}_{V_C}\Big[\frac{|\hat{V}_C\triangle V_C|}{2\vert V_C\vert}\Big],
\end{equation}
%
where $\hat{V}_C\triangle V_C := (V \setminus \hat{V}_C)\cap V_C)\cup( \hat{V}_C \cap (V \setminus V_C))$ denotes the symmetric difference between $\hat{V}_C$ and $V_C$, and $\mathbb{E}_{V_C}$ denotes the expected value given the knowledge of the set $V_C$. Note that $|\hat{V}_C\triangle V_C|\leq 2\vert V_C\vert$ when we assume that the community size is known and $\hat{V}_C$ outputs exactly $k$ vertices, so that in that case $R_{id}\in[0,1]$.
We say that a method achieves \textit{exact recovery} when $R_{id}( \hat{V}_C)\to 0$, and \textit{partial recovery} when $R_{id}( \hat{V}_C)\to c$ for $c\in(0,1)$. In other words, a test achieves partial recovery when it identifies a positive proportion of the vertices in the community. To obtain an estimator for the set of geometric vertices, we construct a test statistic $T: V \to \{0,1\}$ such that $T(i)=1$ if node $i\in V$ is estimated to be in the community, and $T(i)=0$ otherwise.

Low-weight vertices in a GIRG have degree zero with positive probability, and zero-degree type-A and type-B vertices cannot be identified. This strongly suggests that in our setting, where $O(n)$ vertices have weight of order $O(1)$, exact recovery cannot be achieved. In fact, even partial recovery is difficult, because low-weight vertices are a non-vanishing fraction of all the vertices. We therefore focus on achieving partial recovery among the graph induced by all high-weighted vertices. 

For the purpose of identification, we propose the following test statistic $T: V \to \{0,1\}$:
\begin{align}%
T(a) := \mathbbm{1}_{\{ \Wi(a) > n/(w_a \sqrt{\log n})\}}, \label{eq:idteststatistic}
\end{align}%
where
\begin{equation}%
\Wi(a) := \frac{n}{w_a^2}\sum_{b,c \in V}\frac{1}{w_b w_c}\ind{\{a,b,c\} = \triangle}.
\end{equation}%
Next, we show that this leads to vanishing type-I and type-II errors when $w_a = w_a(n) \to\infty$.
\begin{thm}\label{thm:identify}
The test $T(a)$ achieves exact recovery among the set of all vertices $a$ with weight $w_a \gg \log(n)$. Formally, setting
\begin{align}
    \hat{V}_C := \{a\in V : T(a) = 1\},
\end{align}
we get
\begin{align}
   \lim_{n\to\infty}\mathbb E_{V_C}\left[\frac{\vert \hat{V}_C^{t_n} \triangle V_C^{t_n}\vert}{2\vert V_C^{t_n}\vert} \mid \vert V_C^{t_n}\vert >0\right] = 0,
\end{align}
as long as 
\begin{equation}\label{eq:tncondthm}
     (n\log(n)/k)^{1/\tau} \ll t_n \ll k^{1/(\tau-1)}.
\end{equation}
where, for any $U\subseteq V$,
\begin{align}
    U^h:= U \cap \{a\in V : w_a \geq h \}.
\end{align}

\end{thm}

Note that while Theorem~\ref{thm:identify} uses knowledge of $k$ in the threshold for the weights on which the risk tends to 0, this threshold can be avoided by choosing $t_n\gg (n\log(n))^{1/\tau}$, as $k\geq 1$. In this case, one only needs to know whether the upper limit also holds, that is, whether $k\gg (n\log(n))^{(\tau-1)/\tau}$, and one only needs a sufficiently large lower bound on $k$.

\subsection{Estimating the botnet size}
Next we tackle a different issue, namely inferring the size of the planted geometric community under $H_1$. We will make crucial use of the fact that the estimator for the community introduced above identifies high-degree vertices exactly in the $n\to\infty$ limit by Theorem~\ref{thm:identify}. 

Let $X_{(1)},X_{(2)},\dots$ denote the order statistics of the weights of the vertices of the geometric part that are identified by Theorem~\ref{thm:identify}. That is, $X_{(1)}$ is the vertex of the geometric part with the highest degree.   
Thus, we take the $m$ highest-weight vertices that are identified by the node-based test as being part of the GIRG as input. Denote the weights of these vertices by $X_{(1)},X_{(2)}$, and so on.
We propose as an estimator of the community size $k$ the following:
\begin{equation}
    \hat{k}_m := m X_{(m)}^{\tau-1}.\label{eq:kestimators}
\end{equation}

\begin{thm}\label{thm:kestimate}
Assume that $k\gg (n\log(n))^{(\tau-1)/(2\tau-1)}$ and $m\in\mathbb N$. Then, as $n\to\infty$, 
\begin{equation}
    \hat{k}_m/k\plim 1.
\end{equation}
\end{thm}

\begin{proof}[Proof of Theorem~\ref{thm:kestimate}]
Let $X_{(1)},X_{(2)},\dots$ denote the order statistics of the degrees (weights) of the vertices of the GIRG part. By~\cite[Eq.~(4.17)]{Resnick2007}, as $k\to\infty$,
\begin{equation}
    \frac{X_{(s)}}{\Big(\frac{k}{s}\Big)^{1/(\tau-1)}}\plim 1.
\end{equation}
Now the $m$ type-B vertices with highest weight can be identified correctly with probability tending to one as long as there exists some $t_n$ satisfying~\eqref{eq:tncondthm}. Such a sequence exists as long as
\begin{equation}
    k\gg (n\log(n))^{(\tau-1)/(2\tau-1)}.
\end{equation}
Then,
\begin{equation}
    \hat{k}=X_{(m)}^{\tau-1}m/k\plim 1.
\end{equation}
\end{proof}

\subsection{Numerical results}
We present here some numerical experiment to illustrate the finite-sample performance of our tests.

\begin{figure}
     \centering
     \includegraphics[width=\textwidth]{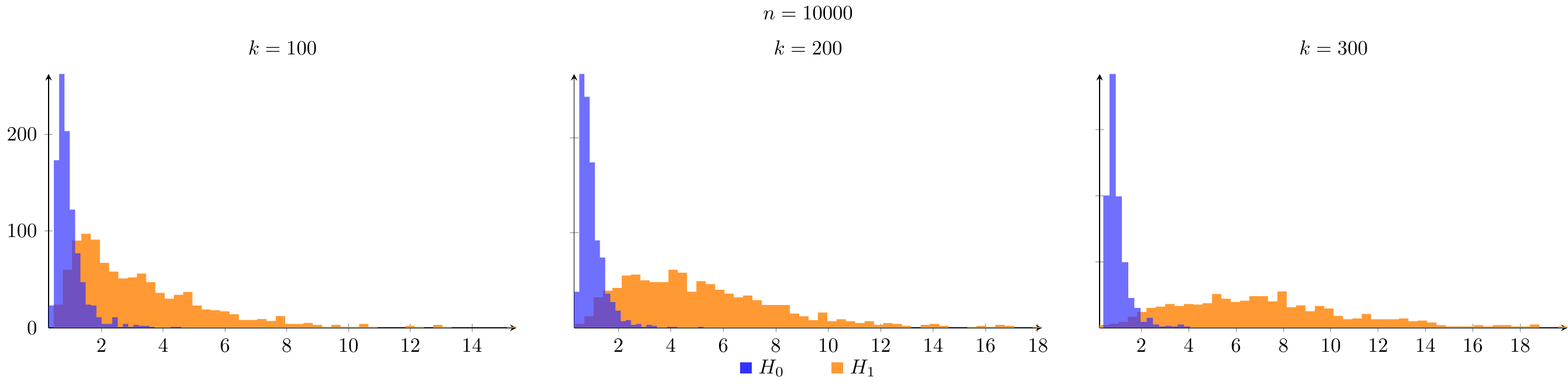}
     \caption{Histogram with the value of $W$ of $10^4$ sample graphs generated under $H_0$ (blue) and under $H_1$ (orange). The parameters chosen for this simulation are $\tau = 2.5, C=1, w_0=1, d=2, \gamma=5$.}
        \label{fig:detection}
\end{figure}

In Figure \ref{fig:detection} we compare the histogram of $W$ from \eqref{eq:W} evaluated over multiple samples of the null and alternative model. In both cases, the models are generated on $n=10^4$ vertices, the size of the geometric community varies between $k=100,200,300$. As we can see, under $H_0$, $W$ is highly concentrated around its expected value. On the other hand, under $H_1$ the typical value of $W$ is larger and increases as the size of the geometric community grows. This is consistent with Proposition \ref{prop:WIRG} and Proposition \ref{prop:WGIRG}, and shows that the weighted triangles test $\psi_W$ can work with high accuracy also for finite samples. Furthermore, the larger $k$ is, the better the separation between $H_0$ and $H_1$ in terms of $W$ is.\\

\begin{figure}
    \centering
    \includegraphics[width=0.6\textwidth]{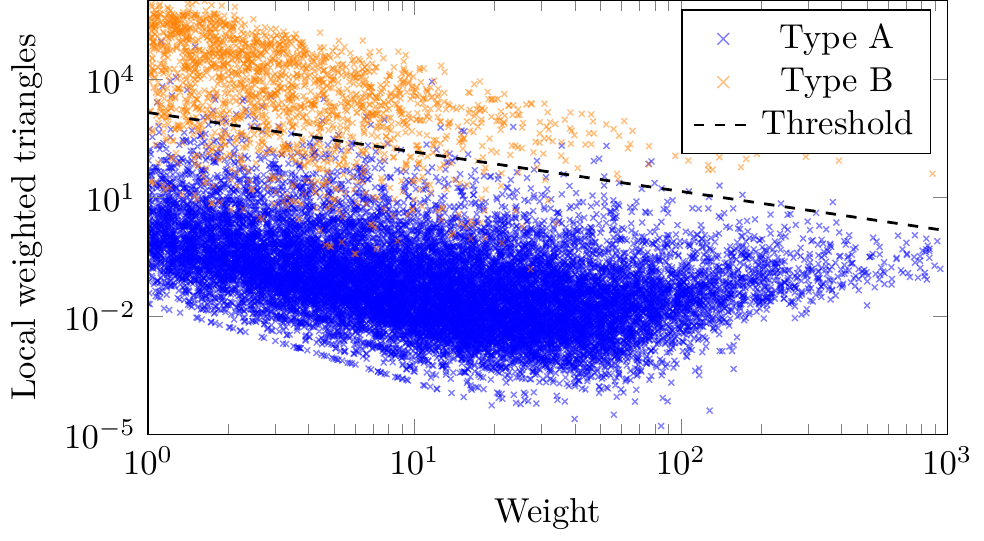}
    \caption{Identification of geometric vertices, under $H_1$. The size of the graph and of the geometric community are respectively $n=10^6$ and $k=10^4$. Dots represent the local weighted triangles statistic for all type-A and type-B vertices, against their weight.}
    \label{fig:identification}
\end{figure}

Next, Figure \ref{fig:identification} illustrates the performance of our identification test. Here, a single instance of $H_1$ has been sampled and the value of the local weighted triangles statistic $W(a)$ for each vertex $a$ is computed. Plotting $W(a)$ against the weights $w_a$, we can notice that the clouds of coordinates $(w_a, W(a))$ of type-B vertices separate from the cloud formed by type-A vertices. The two distinct regions can be easily distinguished in log-log scale. The dotted line in Figure \ref{fig:identification} is the curve $y = \mathcal{C} n/(x \sqrt{\log n})$, where here $\mathcal{C}$ is a constant value tuned on the parameters of the model. According to Theorem \ref{thm:identify}, all but a small fraction of type-B vertices with large weights lie above the dotted line, whereas the type-A vertices are located below. Our simulations confirm this.
We also observe that a large proportion of the vertices with low weights is correctly identified. This suggests that partial recovery could still be achieved for the entire community, even though Theorem~\ref{thm:identify} only works for high-degree vertices, by ignoring the vertices whose local weighted triangles statistic equals zero.\\

\begin{figure}
    \centering
    \includegraphics[width=1\textwidth]{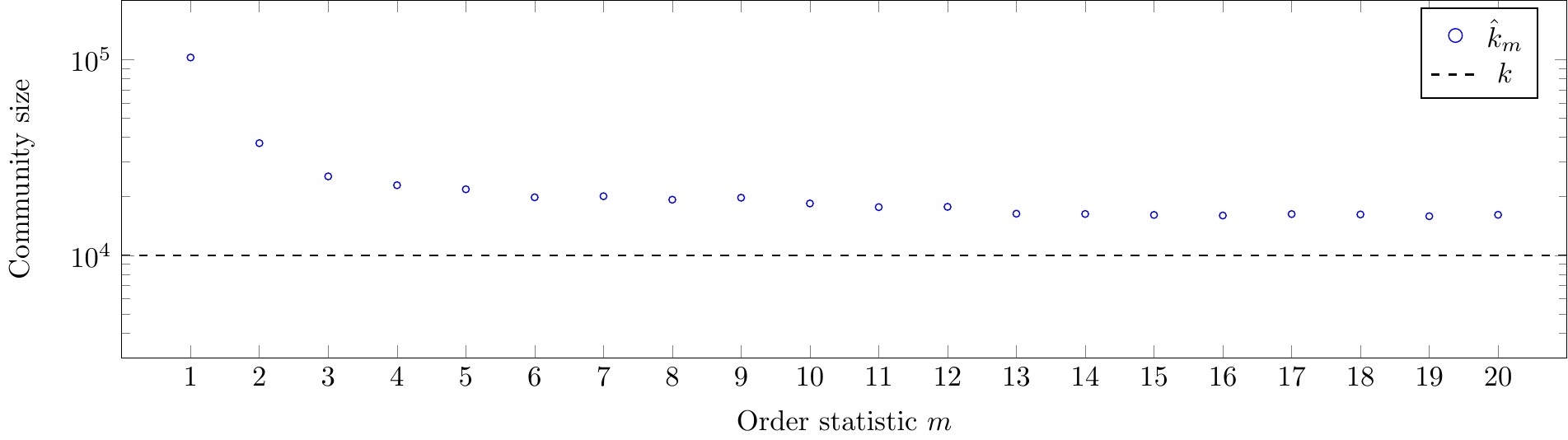}
    \caption{Estimate of the geometric community size. Blue dots are the average over 15 simulations of $\hat{k}_m$. The red line is the value $k$, the true size of the geometric community. The parameters chosen for the model here are: $n=10^6$, $k=10^4$, $\tau=2.5$, $d=1$, $\gamma = 5$, $C=1$, $w_0=1$}
    \label{fig:estimate}
\end{figure}

Lastly, in Figure \ref{fig:estimate} we show a practical application of Theorem \ref{thm:kestimate}. Assume we are given a sample of the model $H_1$. Using the test $T(a)$ defined in \eqref{eq:idteststatistic} and following the procedure mentioned above, the vertices above can be classified either as type-A or type-B. Then, we take the $M$ upper order statistics $(X_{m})_{m =1,...,M}$ for the weights of the vertices classified as type-B, where $M = 20$ in our case. Finally, we compute the community size estimators $(\hat{k}_m)_{m=1,...,M}$, from the definition in \eqref{eq:kestimators}.  Iterating this procedure over multiple independent samples of $H_1$ (keeping all parameters fixed), it is possible to take an average of the size estimators $(\hat{k}_m)_{m=1,...,M}$ for $k$, to check their performance. The blue dots in Figure \ref{fig:estimate} are the average values of the size estimators, obtained from 15 simulations of the model $H_1$ with $n=10^6,k=10^4$. From Theorem \ref{thm:kestimate}, we know that the estimators $\hat{k}_m$ converge in probability to the real size of the geometric community $k$. The convergence cannot be observed in Figure \ref{fig:estimate}, where the size of the graph is fixed ($n = 10^6$). Nonetheless, the figure shows that the estimators $(\hat{k}_m)_m \geq 1$ are able to capture quite nicely the size of the real geometric community $V_C$, even for moderate values of $m$. 

\section{Discussion}\label{sec:discussion}
\paragraph{Computational complexity.} As our test is a triangle-based test, it only requires a triangle enumeration for all vertices. This can be done in $O(n^{1+\frac{1}{\tau}})$~\cite{latapy2008} or $O(n^{3/2})$~\cite{latapy2007practical} time or in time $O(n\log(n))$ for a good approximation~\cite{becchetti2010}, providing an extremely efficient method to detect and identify geometry. 

\paragraph{Iterative procedure.} Our identification procedure identifies high-weight vertices correctly with high probability. We believe that this identification can serve as a starting point to also identify the lower-weight vertices with non-trivial probability. Given the identification of the high-weight neighbors of a low-weight vertex, one can compute the likelihood of this vertex being part of the geometric structure or not, and identify it based on the highest likelihood. After this, one can again assess this likelihood with these updated identifications, and iteratively improve the estimated geometric subset. Such a procedure has been proven to work for the stochastic block model~\cite{yun2014accurate}, and using such procedures in this setting as well seems promising.

\paragraph{Improving  $\hat{k}$.}
While Theorem~\ref{thm:kestimate} shows that our estimator of $k$ is unbiased in the large-network limit, for finite values of $n$, $\hat{k}$ overestimates $k$, as shown in Figure~\ref{fig:estimate}. This is because the test is based on the identified geometric vertices of Theorem~\ref{thm:identify}. In the large-network limit, these are all identified correctly. However, for finite $n$, some vertices may be misclassified as geometric or non-geometric. Since $k$ is small compared to $n$, most misclassifications are non-geometric vertices which are misclassified as geometric. These misidentified vertices therefore make the inferred order statistics of the geometric vertices higher, leading to an overestimation of $k$. Improving this estimate for finite $n$ is therefore an interesting line of further research. For example, one could use information of the expected number of misclassified vertices to improve the test. 

\paragraph{Rescaling the box sizes.}
In our model, we rescaled the GIRG connection probability \eqref{eq:corr_GIRG_connection_rule} with the size of the community $k$. This is equivalent to sampling the locations of the vertices in the community in a shrinking torus, and not rescaling the connection probability. Another natural choice is rescaling the connection probability within the community as
\begin{equation}\label{eq:connection_prob_sparse_alternative}
p_{ij}= \min\left(\frac{w_i w_j}{\mu n ||x_i - x_j||^d}, 1\right)^\gamma.
\end{equation}
However, this connection probability leads to a \textit{sparse community}, where most community members will be disconnected from other community members. 
As the name suggests, in the sparse community scenario the average number of connections between a given vertex in the community and other community vertices decreases roughly as $k/n$. Because of this, we believe that our assumption of a \emph{localized community} hypothesis of \eqref{eq:corr_GIRG_connection_rule} as the more realistic scenario. Still, our detection methods also apply to the sparse community setting as long as $k\gg\sqrt{n}$. The thresholds for identifying such a sparse community are unknown however, and would be an interesting point for further research to investigate the theoretical limits of our methods.
    
\section{Detection: proofs}\label{sec:detect}

In this section we will find an upper bound and a lower bound for the expected weighted triangles in the graph $\Exp{W}$, under the hypothesis $H_0$ and $H_1$. Furthermore, we will derive a bound on the variance of the weighted triangles, which will lead to the proof of Theorem \ref{thm:detect} by applying Chebyshev's theorem. Before computing $\E{W}$ and $\Var(W)$, we state and prove some useful Lemmas.

\begin{lemma}\label{lemma:sumweights}
Let $X$ be a Pareto power-law distribution, with distribution function $F(x) = 1-C x^{1 - \tau}$ for all $x \geq w_0$. Let $\{w_i\}_{i \in V}$ be a weight sequence satisfying Assumption 1. Then,
\begin{enumerate}[(i)]
    \item The sum of the weights is  
    \begin{equation*}
        \sum_{i \in V} w_i = \mu n (1 + o(1))
    \end{equation*}
    where $\mu = \mathbb{E}[X]$.
    \item The sum of the inverse of the weights is 
    \begin{equation*}
        \sum_{i \in V} w_i = \nu n (1 + o(1))
    \end{equation*}
    where $\nu = \mathbb{E}[X^{-1}]$.
\end{enumerate}
\end{lemma}

\begin{proof}
Let $U$ be a uniform distribution over the vertex set $V$. We observe that $F_n$ is the distribution function of a uniformly chosen vertex $w_U$. Indeed,
\begin{equation*}
    \mathbb{P}(w_U \leq x) = \sum_{i \in V} \mathbb{P}(U = i) \mathbbm{1}_{\{w_i \leq x\}} = \sum_{i \in V} \frac{1}{n} \mathbbm{1}_{\{w_i \leq x\}} = F_n.
\end{equation*}
Moreover, we have that
\begin{align*}
    & \mathbb{E}[w_U] = \frac{1}{n} \sum_{i \in V} w_i, \\
    & \mathbb{E}[w_U^{-1}] =  \frac{1}{n} \sum_{i \in V} \frac{1}{w_i}.
\end{align*}
The expected value of $w_U$ can be computed through the cumulative distribution function
\begin{equation}
    \mathbb{E}[w_U] = \int_{\mathbb{R}} x \; \dd F_n(x) = \int_{\mathbb{R}} x (1 + o(1))\; \dd F(x) = \mathbb{E}[X] (1 + o(1)),
\end{equation}
with
\begin{equation}\label{eq:wumin1}
    \mu = \mathbb{E}[X] = \int_{w_0}^{\infty} x \cdot C(\tau -1) x^{-\tau} \; \dd x = C\frac{\tau-1}{\tau-2}w_0^{2-\tau}.
\end{equation}
Similarly, we have $\mathbb{E}[w_U^{-1}] = \mathbb{E}[X] (1 +o(1))$, with 
\begin{equation}\label{eq:invdetailed}
    \nu = \mathbb{E}[X^{-1}] = \int_{w_0}^{\infty} x^{-1} \cdot C(\tau -1) x^{-\tau} \; \dd x = C\frac{\tau-1}{\tau}w_0^{-\tau}.
\end{equation}
\end{proof}

\textit{Remark.} The constant $C$ in the definition of the Pareto random variable $X$ is uniquely determined if we set the minimum weight to be $w_0$. Indeed, $F$ is a proper cumulative distribution function, only when $\int_{w_0}^{\infty} \dd F(x) = 1$. This is the case, when $C = w_0^{\tau-1}$.

\begin{lemma}\label{lemma:sumweightssqr}
\begin{equation*}
    \sum_{i \in V} w_i \mathbbm{1}_{\{w_i \leq \sqrt{n}\}} = \mu n (1 + o(1)).
\end{equation*}
\end{lemma}
\begin{proof}
Following the same proof as in Lemma \ref{lemma:sumweights}, we have that
\begin{align}
    \frac{1}{n}\sum_{i \in V} w_i \mathbbm{1}_{\{w_i \leq \sqrt{n}\}} & =  \int_{w_0}^{\sqrt{n}} x C (\tau-1)x^{-\tau} (1+o(1)) \dd x = [\mu - O(n^{1-\tau/2})](1 +o(1))\nonumber\\
    & = \mu (1 + o(1)). 
\end{align}
\end{proof}


\subsection{Weighted triangles under null hypothesis}

We are now ready to calculate the expectation and variance of the weighted triangles under $H_0$.

\begin{prop}\label{prop:WIRG}
Under $H_0$, the expected value and the variance of $W$ are
\begin{equation}\label{eq:WunderH0quenched}
\begin{split}
    \Expn{W} &= 1+o(1) \\
    \Varn(W) &\leq \frac{1}{\mu^3} (1+o(1)).
\end{split}
\end{equation}
\end{prop}

\begin{proof}
\textbf{Part A: expectation}
\begin{equation}
    \begin{split}
        \Expn{W} 
        &= \sum_{a,b,c \in V} \Expn{\frac{\mathbbm{1}_{\{\{a,b,c\} = \triangle\}}}{w_a w_b w_c}}\\
        &= \sum_{a,b,c \in V} \frac{p_{ab}p_{bc}p_{ac}}{w_a w_b w_c} \label{eq:expWgivenH0}
    \end{split}
\end{equation}
Since $p_{ab} \leq \frac{w_a w_b}{\mu n}$ for any $a,b \in V$,
\begin{equation}
\begin{split}
    \Expn{W} &\leq \left(\frac{1}{\mu n} \right)^3 \sum_{a,b,c \in V} w_a w_b w_c \\
    &= \left(\frac{1}{\mu n} \right)^3 \left(\sum_{i \in V} w_i\right)^3\\
    &= \left(\frac{1}{\mu n} \right)^3 (\mu n)^3 (1 +o(1)) \\
    &= 1+o(1),
\end{split}    
\end{equation}
where $\sum_{i} w_i = \mu n(1 + o(1))$, from Lemma \ref{lemma:sumweights}(i).
To obtain a lower bound, we restrict the sum in \eqref{eq:expWgivenH0} to those $a,b,c$ such that $w_a,w_b,w_c\leq \sqrt{\mu n}$. Then, $\frac{w_a w_b}{\mu n},\frac{w_b w_c}{\mu n},\frac{w_a w_c}{\mu n} \leq 1$, and applying Lemma \ref{lemma:sumweightssqr},
\begin{equation}
\begin{split}
    \Expn{W} 
    &\geq \sum_{a : w_a \leq \sqrt{\mu n}} \; \sum_{b : w_b \leq \sqrt{\mu n}} \; \sum_{c : w_c \leq \sqrt{\mu n}} \frac{p_{ab}p_{bc}p_{ac}}{w_a w_b w_c} \\
    &= \left(\frac{1}{\mu n}\right)^3 \sum_{a : w_a \leq \sqrt{\mu n}} \; \sum_{b : w_b \leq \sqrt{\mu n}} \; \sum_{c : w_c \leq \sqrt{\mu n}} w_a w_b w_c \\
    &= \left(\frac{1}{\mu n}\right)^3 (n \mu )^3(1+o(1)) \\
    &= 1 +o(1),
\end{split}
\end{equation}
where we applied Lemma \ref{lemma:sumweightssqr}.\\

\noindent\textbf{Part B: variance}

\noindent We can rewrite the variance of $W$ as
\begin{equation}
\begin{split}\label{eq:var_h0}
    \Varn(W) &= \Varn \left(\sum_{a,b,c \in V} \frac{\mathbbm{1}_{\{\{a,b,c\} = \triangle\}}}{w_a w_b w_c} \right)\\ 
    &=\sum_{a_1,b_1,c_1 \in V} \sum_{a_2,b_2,c_2 \in V} \frac{\Cov_0\left(\mathbbm{1}_{\{\{a_1,b_1,c_1\} = \triangle\}}, \mathbbm{1}_{\{\{a_2,b_2,c_2\} = \triangle\}}]\right)}{w_{a_1} w_{b_1} w_{c_1} w_{a_2} w_{b_2} w_{c_2}}
\end{split}
\end{equation}
Then, using the bound $\Cov(X_i,X_j) \leq \Exp{X_i X_j}$ for the covariance,

\begin{equation}
\begin{split}
    \Varn(W) &\leq \sum_{a_1,b_1,c_1 \in V} \sum_{a_2,b_2,c_2 \in V} \frac{\Expn{\mathbbm{1}_{\{\{a_1,b_1,c_1\} = \triangle\}} \mathbbm{1}_{\{\{a_2,b_2,c_2\} = \triangle\}}}}{w_{a_1} w_{b_1} w_{c_1} w_{a_2} w_{b_2} w_{c_2}}\\
    &=\sum_{a_1,b_1,c_1 \in V} \; \sum_{a_2,b_2,c_2 \in V} \frac{P_0(a_1 \leftrightarrow b_1,b_1 \leftrightarrow c_1,c_1 \leftrightarrow a_1,a_2 \leftrightarrow b_2,b_2 \leftrightarrow c_2,c_2 \leftrightarrow a_2)}{w_{a_1} w_{b_1} w_{c_1}w_{a_2} w_{b_2} w_{c_2}}
    \label{eq:varboundWgivenH0}
\end{split}
\end{equation}
There are now different cases for the intersection of the six vertices $a_1,b_1,c_1,a_2,b_2,c_2$:
\begin{itemize}
    \item If $\{a_1,b_1,c_1\}$ and $\{a_2,b_2,c_2\}$ do not intersect or intersect in one vertex only, then their contribution to the variance in~\eqref{eq:var_h0} is zero. Indeed, in this case the two Bernoulli random variables $\mathbbm{1}_{\{\{a_1,b_1,c_1\} = \triangle\}},\mathbbm{1}_{\{\{a_2,b_2,c_2\} = \triangle\}}$ are uncorrelated. 
    
    \item If $\{a_1,b_1,c_1\}$ and $\{a_2,b_2,c_2\}$ intersect in 2 vertices, without loss of generality, and introducing a combinatorial factor, we may assume $a_1= a_2\equiv a, b_1 = b_2 \equiv b$.
    Then, the numerator in the sum of \eqref{eq:varboundWgivenH0} is bounded by
    \begin{equation}
        \begin{split}
            P_0(a \leftrightarrow b, a \leftrightarrow c_1, b \leftrightarrow c_1, a \leftrightarrow c_2, b \leftrightarrow c_2) \leq \frac{w_a^3 w_b^3 w_{c_1}^2 w_{c_2}^2}{(\mu n)^5},
        \end{split}
    \end{equation}
    using the fact that $p_{ij} \leq w_i w_j / \mu n$, for all $i,j$. Then the contribution to the variance from such vertices is bounded by 
    \begin{equation}
        \frac{1}{(\mu n)^5}\sum_{a,b,c_1,c_2 \in V} \frac{w_a^3 w_b^3 w_{c_1}^2 w_{c_2}^2} {w_{a}^2 w_{b}^2 w_{c_1} w_{c_2}} = O(n^{-1})
    \end{equation}
    where the equality follows from Lemma \ref{lemma:sumweights}.
    
    \item If $\{a_1,b_1,c_1\}$ and $\{a_2,b_2,c_2\}$ intersect in all 3 vertices, then without loss of generality, and up to a combinatorial factor, we may assume that $a_1= a_2 \equiv a, b_1 = b_2 \equiv b, c_1 = c_2 \equiv c$.
    In this case the numerator in the sum of \eqref{eq:varboundWgivenH0} is bounded by
    \begin{equation}
        P_0(a \leftrightarrow b, a \leftrightarrow c, b \leftrightarrow c) \leq \frac{w_a^2 w_b^2 w_{c}^2}{(\mu n)^3}.
    \end{equation}
    Therefore, the contribution to the variance from such vertices is bounded by
    \begin{equation}
        \frac{1}{(\mu n)^3}\sum_{a,b,c \in V} \frac{w_a^2 w_b^2 w_{c}^2} {w_{a}^2 w_{b}^2 w_{c}^2} = O(1).
    \end{equation}
\end{itemize}

Summing up:
\begin{equation}
    \Expn{W} = \frac{1}{6} (1+o(1)),
\end{equation}
and 
\begin{equation}
    \Varn(W) = O(1).
\end{equation}
\end{proof}

\subsection{Weighted triangles under alternative hypothesis}
Next, we bound the expected weighted triangles and their variance under the alternative hypothesis $H_1$. Remember that, under $H_1$, there exist two different types of vertices, type-A and type-B, which are part respectively of the non-geometric part of the graph or the geometric community.
Triangles in the alternative model are then of 4 different types: 
\begin{itemize}
    \item those between type-A vertices, denoted by $\triangle_1$ 
    \item those with two type-A vertices, one type-B vertex, denoted by $\triangle_2$
    \item those with one type-A vertex, two type-B vertices, denoted by $\triangle_3$
    \item those between type-B vertices, denoted by $\triangle_4$
\end{itemize}

\begin{prop}\label{prop:WGIRG}
Assume $k \equiv k(n) \to \infty$, as $n \to \infty$. Under the hypothesis $H_1$, the expected value and the variance of $W$ are
\begin{equation}\label{eq:WunderH0annealed}
    \begin{split}
        \Expa{W} &= \Omega(k) \\
        \Vara(W) &= O(k) 
    \end{split}
\end{equation}
\end{prop}

\begin{proof}

\textbf{Part A: expectation}

We can split the expected value of $W$ into the sum over all possible types of triangles:
\begin{equation}
    \Expa{W} = \Expa{\sum_{\alpha = 1}^4 \sum_{\{a,b,c\} \subset V} \frac{\ind{\{a,b,c\} \in \triangle_{\alpha}}}{w_a w_b w_c} } = \sum_{\alpha = 1}^4 \sum_{\{a,b,c\} \subset V} \frac{P_1(\{a,b,c\} \in \triangle_{\alpha})}{w_a w_b w_c}
\end{equation}

Then, we lower bound the expected value of $W$, with the contribution from $\triangle_4$:
\begin{equation} \label{eq:expWgivenH1}
\begin{split}
     \Expa{W} \geq \sum_{\{a,b,c\} \subset V} \frac{P_1(\{a,b,c\} \in \triangle_4)}{w_a w_b w_c} = \sum_{\{a,b,c\} \subset V_C} \frac{P_1(a \leftrightarrow b, b \leftrightarrow c, c \leftrightarrow a)}{w_a w_b w_c}
\end{split}
\end{equation}
Let $\{a,b,c\} \subset V_C$ be fixed now. Since $a,b,c$ are type-B vertices,
\begin{equation}\label{eq:pt4int}
\begin{split}
     & P_1(a \leftrightarrow b, b \leftrightarrow c, c \leftrightarrow a) = \mathbb{E}_{\overline{x}} \left[ p(w_a,w_b,x_a,x_b) p(w_b,w_c,x_b,x_c) p(w_a,w_c,x_a,x_c) \right] \\
     &= \int_{(\mathbb{T}^d)^3} \left(\frac{w_a w_b}{\mu k ||x_a - x_b||^d} \wedge 1 \right)^{\gamma} \left(\frac{w_b w_c}{\mu k ||x_b - x_c||^d} \wedge 1 \right)^{\gamma} \left(\frac{w_a w_c}{\mu k ||x_a - x_c||^d} \wedge 1 \right)^{\gamma} \dd \overline{x}
\end{split}
\end{equation}
where $\overline{x} = \{x_a,x_b,x_c\}$ and $\mathbb{E}_{\overline{x}}$ is the expectation over the positions $\overline{x}$. Observe that 
\begin{equation}
    \begin{split}
        ||x_b - x_c|| &\leq ||x_b - x_a|| + ||x_a - x_c|| \\
        &\leq 2 \max\{||x_a - x_b||, ||x_a - x_c||\}.
    \end{split}
\end{equation} 
Thus,
\begin{equation}
    \begin{split}
        \frac{1}{||x_b - x_c||^d} &\geq \frac{1}{2^d \max\{||x_a - x_b||^d, ||x_a - x_c||^d\}} \\&= \min\left\{ \frac{1}{2^d ||x_a - x_b||^d}, \frac{1}{2^d ||x_a - x_c||^d} \right\}.
    \end{split}
\end{equation}
Then, using the substitution $\{x'_a, x'_b, x'_c\} = \{x_a, x_b - x_a, x_c - x_a\}$,~\eqref{eq:pt4int} is lower bounded by
\begin{equation}
    \int_{(\mathbb{T}^d)^2} \left(\frac{w_a w_b}{\mu k ||x'_b||^d} \wedge 1 \right)^{\gamma} \left(\frac{w_a w_c}{\mu k ||x'_c||^d} \wedge 1 \right)^{\gamma} \left(\frac{w_b w_c}{\mu k 2^d||x'_b||^d} \wedge \frac{w_b w_c}{\mu k 2^d||x'_c||^d} \wedge 1 \right)^{\gamma} \dd x'_b \; \dd x'_c  \label{eq:threeintegrands}
\end{equation}
Next, we represent the torus $\mathbb{T}^d$ as the interval $[-\frac{1}{2},\frac{1}{2}]^d$, and consider the cube around 0
\begin{equation}
    \mathcal{C} = \left[-\left(\frac{w_0^2}{\mu k 2^d}\right)^{1/d}, \left(\frac{w_0^2}{\mu k 2^d}\right)^{1/d}\right]^d,
\end{equation}
where $w_0$ is the minimum weight in the model. If $x'_b ,x'_c \in \mathcal{C}$, then all the minima in the integrand of \eqref{eq:threeintegrands} are equal to 1. Therefore, if we restrict the integral domain to the subset $\mathcal{C}^2$, \eqref{eq:threeintegrands} is lower bounded by the volume of $\mathcal{C}^2$. Summing up,
\begin{equation}
        P_1(a \leftrightarrow b, b \leftrightarrow c, c \leftrightarrow a) \geq \int_{\mathcal{C}^2} \dd x'_b \; \dd x'_c = \left(\frac{w_0^2}{\mu k}\right)^2
\end{equation}
Hence,
\begin{equation}
    \begin{split}
        \Expa{W} &\geq \sum_{\{a,b,c\} \subset V_C} \frac{P_1(a \leftrightarrow b, b \leftrightarrow c, c \leftrightarrow a)}{w_a w_b w_c} \\
        &\geq \left(\frac{w_0^2}{\mu k}\right)^2 \sum_{\{a,b,c\} \subset V_C} \frac{1}{ w_a w_b w_c} \\
        &= \left( \frac{w_0^2}{\mu k}\right)^2  \Theta(k^3)\\
        &= \Theta(k),
    \end{split}
\end{equation}
where, in the equality, we applied Lemma \ref{lemma:sumweights}(ii). In conclusion, $\Expa{W} = \Omega(k)$. \\

\noindent \textbf{Part B: variance}

\noindent  We can upper bound the variance of $W$ under $H_1$ as follows:
\begin{align}
    \Vara(W) &= \Vara\left(\sum_{a,b,c \in V} \frac{\mathbbm{1}_{\{\{a,b,c\} = \triangle\}}}{w_a w_b w_c} \right)\nonumber\\
    &= \sum_{a_1,b_1,c_1 \in V} \; \sum_{a_2,b_2,c_2 \in V} {\frac{\Cov_1\left(\mathbbm{1}_{\{\{a_1,b_1,c_1\} = \triangle\}} ,\mathbbm{1}_{\{\{a_2,b_2,c_2\} = \triangle\}}\right)}{w_{a_1} w_{b_1} w_{c_1} w_{a_2} w_{b_2} w_{c_2}}}\nonumber\\
    &\leq \sum_{a_1,b_1,c_1 \in V} \; \sum_{a_2,b_2,c_2 \in V} {\frac{\Expa{\mathbbm{1}_{\{\{a_1,b_1,c_1\} = \triangle\}} \mathbbm{1}_{\{\{a_2,b_2,c_2\} = \triangle\}}}}{w_{a_1} w_{b_1} w_{c_1} w_{a_2} w_{b_2} w_{c_2}}}\nonumber\\
    &= \sum_{a_1,b_1,c_1 \in V} \; \sum_{a_2,b_2,c_2 \in V} \frac{P_1(\{a_1,b_1,c_1\} = \triangle, \{a_2,b_2,c_2\} = \triangle)}{w_{a_1} w_{b_1} w_{c_1}w_{a_2} w_{b_2} w_{c_2}} .\label{eq:varboundWgivenH1}
\end{align}
The terms in the last sum of equation \eqref{eq:varboundWgivenH1} depend on the vertex types of $\{a_1,b_1,c_1\}$ and $\{a_2,b_2,c_2\}$, as the vertex types determine the connection probabilities.

First, consider the case when both $\{a_1, b_1, c_1\}$ and $\{a_2, b_2, c_2\}$ contain at most one type-B vertex. In this case, all connections are non-geometric. Therefore, the contribution to the variance from this combinations of vertices is $O(1)$, following the proof of Proposition \ref{prop:expWiIRG} which bounds the variance of non-geometric triangles. 

Moreover, for the set of vertices such that $\{a_1, b_1, c_1\}$ and $\{a_2, b_2, c_2\}$ do not intersect, the random variables $\mathbbm{1}_{\{\{a_1,b_1,c_1\} = \triangle\}}$ and $\mathbbm{1}_{\{\{a_2,b_2,c_2\} = \triangle\}}$ are independent Bernoulli random variables, and their contribution to the variance is 0. \\

Thus, we only focus on the case when at least one of the sets $\{a_1, b_1, c_1\}$ and $\{a_2, b_2, c_2\}$ contains two or three type-B vertices, and when the two sets intersect in at least one vertex.

\begin{itemize}
    \item If $|\{a_1,b_1,c_1\} \cap \{a_2,b_2,c_2\}| = 1$, without loss of generality, up to a combinatorial factor, $a_1 = a_2 \equiv a$. In this case, at least one of the paths $(c_1,b_1,a,b_2,c_2)$ or $(b_1,c_1,a,c_2,b_2)$ contains at least one geometric edge, that is, it has at least one pair of adjacent type-B vertices. 
    \begin{figure}[h!]
        \centering
        \includegraphics[scale=0.75]{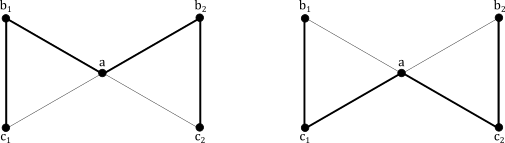}
    \end{figure}
    Suppose that the path $(c_1,b_1,a,b_2,c_2)$ contains $m \geq 1$ geometric edges (the proof is analogous in the other case). Then, the probability in the sum of \eqref{eq:varboundWgivenH1} is bounded by
    \begin{align}
        P_1(\{a,b_1,c_1\} = \triangle, \{a,b_2,c_2\} = \triangle) &\leq P_1(c_1 \leftrightarrow b_1 \leftrightarrow a \leftrightarrow b_2 \leftrightarrow c_2)\nonumber\\
        &\leq \left(\frac{1}{1+C_1}\right)^4 \frac{w_a^2 w_{b_1}^2 w_{b_2}^2 w_{c_1} w_{c_2}}{\mu^4 k^m n^{4-m}},
    \end{align}
    where the last inequality follows from the fact that $m$ of the four edges have geometric connection rule, combined with Lemma \ref{lemma:marginalprob} and Lemma \ref{lemma:geopathbound}. Since there are $m$ geometric edges in the path, at least $m+1$ of the vertices $a,b_1,b_2,c_1,c_2$ need to be type-B. Then, $m+1$ vertices need to be chosen in at most $k$ different ways, the remaining in at most $n$ different ways. Therefore, the contribution in \eqref{eq:varboundWgivenH1} from all possible combinations of vertices in this case is
    \begin{align}
        &\sum_{\substack{a,b_1,b_2,c_1,c_2 \in V \\ m+1 \text{ are type-B} }} \frac{P(\{a,b_1,c_1\} = \triangle, \{a,b_2,c_2\} = \triangle)}{w_{a}^2 w_{b_1} w_{c_1} w_{b_2} w_{c_2}}\nonumber\\
        &\leq \sum_{\substack{a,b_1,b_2,c_1,c_2 \in V \\ m+1 \text{ are type-B} }} O\left(\frac{w_a^2 w_{b_1}^2 w_{b_2}^2 w_{c_1} w_{c_2}}{k^m n^{4-m} w_a^2 w_{b_1} w_{c_1} w_{b_2} w_{c_2}}\right) \nonumber\\
        &= \sum_{\substack{a,b_1,b_2,c_1,c_2 \in V \\ m+1 \text{ are type-B} }} O\left(\frac{w_{b_1} w_{b_2}}{k^m n^{4-m}}\right) \nonumber\\
        &= O(k),
    \end{align}
    where last equality follows from Lemma \ref{lemma:sumweights}.

    Summing over all possible choices for $m = 1,2,3,4$ we end that the contribution to the variance from the the triplets $\{a_1,b_1,c_1\}$, $\{a_2,b_2,c_2\}$ intersecting in one vertex is $O(k)$.
    
    \item If $|\{a_1,b_1,c_1\} \cap \{a_2,b_2,c_2\}| = 2$, without loss of generality $a_1 = a_2 \equiv a$ and $b_1 = b_2 \equiv b$. In this case, at least one of the paths $(c_1,b,a,c_2)$ or $(c_1,a,b,c_2)$ has at least one geometric edge.
    \begin{figure}[h!]
        \centering
        \includegraphics[scale=0.75]{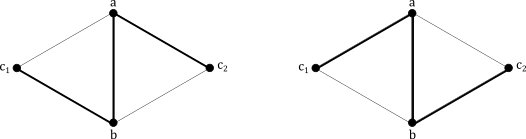}
    \end{figure}
    Assume the path $(c_1,b,a,c_2)$ contains $m \geq 1$ geometric edges. Then, similarly as before, we can bound the probability in \eqref{eq:varboundWgivenH1} by
    \begin{align}
        P_1(\{a,b,c_1\} = \triangle, \{a,b,c_2\} = \triangle) &\leq P_1(c_1 \leftrightarrow b \leftrightarrow a \leftrightarrow c_2)\nonumber\\
        &\leq \left(\frac{1}{1+C_1}\right)^3 \frac{w_a^2 w_{b}^2 w_{c_1} w_{c_2}}{\mu^3 k^m n^{3-m}}.
    \end{align}
    As before, $m+1$ of the vertices need to be type-B. Therefore, the contribution in \eqref{eq:varboundWgivenH1} from all possible combinations of vertices in this case is again
    \begin{align}
        \sum_{\substack{a,b,c_1,c_2 \in V \\ m+1 \text{ are type-B} }} \frac{P_1(\{a,b,c_1\} = \triangle, \{a,b,c_2\} = \triangle)}{w_{a}^2 w_{b}^2 w_{c_1} w_{c_2}} &\leq \sum_{\substack{a,b,c_1,c_2 \in V \\ m+1 \text{ are type-B} }} O\left(\frac{w_a^2 w_{b}^2 w_{c_1} w_{c_2}}{k^m n^{3-m} w_a^2 w_{b}^2 w_{c_1} w_{c_2}}\right) \nonumber\\
        &= \sum_{\substack{a,b,c_1,c_2 \in V \\ m+1 \text{ are type-B} }} O\left(\frac{1}{k^m n^{4-m}}\right) \nonumber\\
        &= O(k).
    \end{align}
    Summing over all possible choices for $m = 1,2,3$ we end that the contribution to the variance from the the triplets $\{a_1,b_1,c_1\}$, $\{a_2,b_2,c_2\}$ intersecting in two vertices is $O(k)$.
    
    \item If $|\{a_1,b_1,c_1\} \cap \{a_2,b_2,c_2\}| = 3$, without loss of generality $a_1 = a_2 \equiv a, b_1 = b_2 \equiv b, c_1 = c_2 \equiv c$. In this case, one of the paths $(a,b,c)$ or $(a,c,b)$ has at least one pair of adjacent type-B vertices. 
    \begin{figure}[h!]
        \centering
        \includegraphics[scale=0.75]{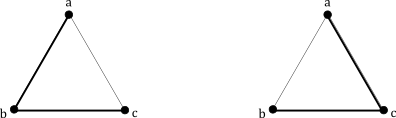}
    \end{figure}
    Assume it is the path $(a,b,c)$ having $m$ geometric edges.
    Following the same steps as before, we bound the probability in \eqref{eq:varboundWgivenH1} by
    \begin{align}
        P_1(\{a,b,c\} = \triangle) &\leq P_1(a \leftrightarrow b \leftrightarrow c)\\
        &\leq \left(\frac{1}{1+C_1}\right)^2 \frac{w_a w_{b}^2 w_c}{\mu^2 k^m n^{2-m}}.
    \end{align}
    Since at least $m+1$ of the three vertices are type-B, then the contribution in \eqref{eq:varboundWgivenH1} from all possible combinations of vertices in this case is once again
    \begin{align}
        \sum_{\substack{a,b,c \in V \\ m+1 \text{ are type-B} }} \frac{P_1(\{a,b,c\} = \triangle)}{w_{a}^2 w_{b}^2 w_{c}^2} &\leq \sum_{\substack{a,b,c \in V \\ m+1 \text{ are type-B} }} O\left(\frac{w_a w_{b}^2 w_c}{k^m n^{2-m} w_a^2 w_{b}^2 w_{c}^2}\right) \nonumber\\
        &= \sum_{\substack{a,b,c \in V \\ m+1 \text{ are type-B} }} O\left(\frac{1}{k^m n^{2-m} w_a w_c}\right) \nonumber\\
        &= O(k).
    \end{align}
    Then, summing over all possible choices for $m = 1,2$ we find that the contribution to the variance from the the triplets $\{a_1,b_1,c_1\}$, $\{a_2,b_2,c_2\}$ intersecting in three vertices is $O(k)$.
\end{itemize}
Summing up all possible cases, we conclude that the variance is $O(k)$.\\
\end{proof}

\begin{proof}[Proof of Theorem~\ref{thm:detect}]
    By Proposition~\ref{prop:WIRG} and Chebyshev's inequality,
    \begin{equation}
        P_0(W(G)>f(n))\leq \frac{\Var_0(W(G))}{f(n)^2} \longrightarrow 0.
    \end{equation}
    Furthermore, by Proposition~\ref{prop:WGIRG}, for $n$ sufficiently large, there exists a constant $M$ such that
    \begin{equation}
        P_1(W(G)<f(n))\leq P_1(W(G)-\mathbb E[W(G)]<f(n)-Mk).
    \end{equation}
    Then, again by Chebyshev's inequality,
    \begin{equation}
        P_1(W(G)<f(n))\leq \frac{\Var_1(W(G))}{(f(n)-Mk)^2}\to 0,
    \end{equation}
    by Proposition~\ref{prop:WGIRG}.
\end{proof}

\section{Localized triangle statistic: proof of Theorem~\ref{thm:identify}}\label{sec:identify}

We will now show that, under the alternative hypothesis $H_1$, a positive fraction of high-degree vertices in the geometric community of the network can be distinguished though the localized triangle statistic.

\subsection{Localized weighted triangle for vertices in the IRG}
\begin{prop}\label{prop:expWiIRG}
Suppose $a$ is a type-A vertex. When $k=o(n)$, 
\begin{equation}
    \Expa{\Wi(a)}\leq\frac{1}{2\mu}(1+o(1)).
\end{equation}
\end{prop}

\begin{proof}
We can write
	\begin{align}
		\Expa{\Wi(a)}& = \frac{n}{2 w_a^2}\sum_{b,c \in V} \Expa{\frac{1}{w_b w_b}\ind{\{a,b,c\}=\triangle}}\nonumber\\
		&=\frac{n}{2 w_a^2}\sum_{b,c \in V} \frac{p_{ab}p_{bc}p_{ac}}{w_b w_c},
	\end{align}
because the edges $\{a,b\},\{b,c\},\{a,c\}$ are independently present in the graph. Since $a$ is type-A, we can bound the connection probabilities $p_{ab} \leq w_a w_b/\mu n$ and $p_{ac} \leq w_a w_c/\mu n$. When at least one of $b$ or $c$ is a type-A vertex, then the connection probability $p_{bc}$ is bounded by $w_b w_c/\mu n$. In this case,
\begin{align}
    \frac{n}{2 w_a^2}\sum_{\substack{b,c \in V \\ b \text{ or } c \text{ type-A}}} \frac{p_{ab}p_{bc}p_{ac}}{w_b w_c} &\leq \frac{n}{2 w_a^2}\sum_{b,c \in V} \frac{w_a^2 w_b w_c}{(\mu n)^3} \nonumber \\
    &= \frac{1}{2 \mu^3 n^2}\sum_{b,c \in V} {w_b w_c} \nonumber \\
    &\leq \frac{1}{2 \mu^3 n^2} (\mu n)^2 (1 + o(1)) = \frac{1}{2 \mu}(1 + o(1)).
\end{align}
On the other hand, when both $b$ and $c$ are type-B vertices, the connection probability $p_{bc}$ is $O(w_b w_c/k)$. Then,
\begin{align}
    \frac{n}{2 w_a^2}\sum_{\substack{b,c \in V \\ b \text{ and } c \text{ type-B}}} \frac{p_{ab}p_{bc}p_{ac}}{w_b w_c} &\leq \frac{n}{2 w_a^2}\sum_{b,c \in V \setminus [n-k]} O\left(\frac{w_a^2 w_b w_c}{\mu^3 n^2 k}\right) \nonumber \\
    &= O\left(\frac{k}{2 \mu n}\right) = o(1).
\end{align}
\end{proof}

\subsection{Node-based statistics for node in the geometric community}
We now consider the node-based triangle statistic for a node in the geometric community. 

\begin{prop}\label{prop:expWiGIRG}
Suppose $a$ is a type-B vertex. When $k=o(n)$ and $k\to\infty$ as $n\to\infty$,
\begin{equation}
\Expa{\Wi(a)} = \Omega(n/w_a)
\end{equation}
\end{prop}

\begin{proof}
Assume $w_a \gg \log(n)$. We bound the node-based statistic by only the triangles with vertex $a$ where all other vertices $b,c$ of the triangle are also type-B.
\begin{equation}
\begin{split}\label{eq:W(a)firstbound}
    \Expa{\Wi(a)} &\geq \frac{n}{2 w_a^2}\sum_{b,c \in V_C} \Expa{\frac{1}{w_b w_c}\ind{\{a,b,c\}=\triangle}}\\
	&=\frac{n}{2 w_a^2}\sum_{b,c \in V_C} \frac{P_1(\{a,b,c\}=\triangle)}{w_b w_c}\\
	&=\frac{n}{2 w_a^2}\sum_{b,c \in V_C} \frac{\mathbb{E}_{x}\left[ p_{ab} p_{bc} p_{ac}\right]}{w_b w_c}
\end{split}
\end{equation}
where $p_{ij}$ for any $i,j$ type-B vertices is the connection probability defined in \eqref{eq:corr_GIRG_connection_rule}. Observe that for any $i,j$ type-B vertices we can lower bound $p_{ij}$ by
\begin{equation}
    p_{ij} \geq \Tilde{p}_{ij} := 
    \begin{cases}
    \frac{1}{1+C_1}, &\text{ if $||x_i - x_j||^d \leq w_i w_j/\mu k$} \\
    0, &\text{ otherwise.}
    \end{cases}
\end{equation}
Therefore, from \eqref{eq:W(a)firstbound}
\begin{equation}
    \Expa{\Wi(a)} \geq \frac{n}{2 w_a^2(1+C_1)^3}\sum_{b,c \in V_C} \frac{\mathbb{P}_1\left(||x_a-x_b||^d \leq \frac{w_a w_b}{\mu k}, ||x_b-x_c||^d \leq \frac{w_b w_c}{\mu k}, ||x_a-x_c||^d \leq \frac{w_a w_c}{\mu k} \right)}{w_b w_c}. \label{eq:lowerprobtriangle}
\end{equation}

\begin{figure}[tbp]
    \centering
    \includegraphics[scale=0.5]{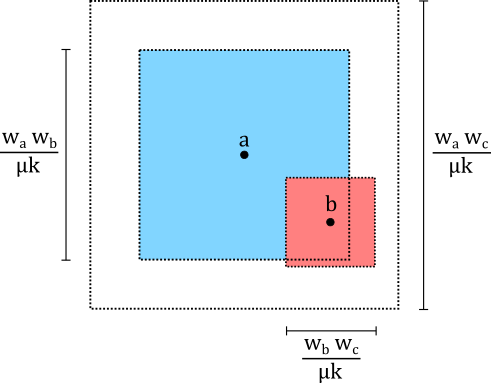}
    \caption{Visual description of the probability in Equation \ref{eq:lowerprobtriangle}, when $d=2$.}
    \label{fig:probvisualproof}
\end{figure}
Assuming $w_b < w_c \leq \log(n)$, the probability in the latter sum is proportional to $\frac{w_a w_b}{\mu k}\frac{w_b w_c}{\mu k} = \frac{w_a w_b^2 w_c}{(\mu k)^2}$. For a visual insight of the latter statement, look at the picture in Figure \ref{fig:probvisualproof}: select the position of $a$; $a$ connects to $b$ when the position of $b$ is inside the blue square; next, $c$ connect to both $a$ and $b$ when its position is in the intersection of the red and white square, which has volume proportional to the (asymptotically smaller) red square.

Then, for some $K>0$,
\begin{equation}
\begin{split}
    \Expa{\Wi(a)} &\geq \frac{nK}{2 w_a^2(1+C_1)^3}\sum_{\substack{b,c \in V_C\\ w_b<w_c \leq \log(n)}} \frac{1}{w_b w_c}\frac{w_a w_b^2 w_c}{(\mu k)^2} \\
    &\geq \frac{nK}{2 w_a(1+C_1)^3} \sum_{\substack{b,c \in V_C\\ w_b<w_c \leq \log(n)}} \frac{w_0}{(\mu k)^2} \\
    &= \frac{w_0}{(1+C_1)^3\mu^2} \frac{nK}{w_a} (1 + o(1)).
\end{split}
\end{equation}
\end{proof}

\subsection{Variance of the localized statistic}

\begin{prop}\label{prop:varianceWi}
Suppose $k=o(n)$, and $k\to\infty$ as $n\to\infty$. For any $a$,
\begin{equation}
    \Var_1(\Wi(a)) = O(n^2/w_a^3).
\end{equation}
\end{prop}

\begin{proof}
The variance of the localized statistic computed on any vertex $a$ is
\begin{align}
    \Var_1(\Wi(a)) &= \frac{n^2}{w_a^4}\sum_{b_1,c_1,b_2,c_2}\Cov_1 \left( \frac{\mathbbm{1}_{\triangle_{a,b_1,c_1}}}{w_{b_1} w_{c_1}}, \frac{\mathbbm{1}_{\triangle_{a,b_2,c_2}}}{w_{b_2} w_{c_2}} \right) \nonumber\\
    &= \frac{n^2}{w_a^4}\sum_{b_1,c_1,b_2,c_2}\frac{P_1(\triangle_{a,b_1,c_1},\triangle_{a,b_2,c_2})-  P(\triangle_{a,b_1,c_1})P(\triangle_{a,b_2,c_2})}{w_{b_1}w_{b_2}w_{c_1}w_{c_2}} \nonumber\\
    &\leq \frac{n^2}{w_a^4}\sum_{b_1,c_1,b_2,c_2}\frac{P_1(\triangle_{a,b_1,c_1},\triangle_{a,b_2,c_2})}{w_{b_1}w_{b_2}w_{c_1}w_{c_2}}.\label{eq:locvarbound}
\end{align}
From the second line of \eqref{eq:locvarbound}, we see that the contribution to the variance from disjoint combinations of vertices $\{b_1,c_1\}$,$\{b_2,c_2\}$ is always zero, as the random variables $\triangle_{a,b_1,c_1},\triangle_{a,b_2,c_2}$ are independent.

Thus, we only focus on the summation terms of the upper bound in \eqref{eq:locvarbound} when $\{b_1,c_1\}$,$\{b_2,c_2\}$ intersect.

\begin{itemize}
    \item If $|\{b_1,c_1\} \cap \{b_2,c_2\}| = 1$, w.l.o.g. (and up to a combinatorial factor) we may assume $b_1 = b_2 \equiv b$.
    Observe that the numerator in the summation term of \eqref{eq:locvarbound} can be further upper bounded by
    \begin{equation}
        P_1(\triangle_{a,b,c_1},\triangle_{a,b,c_2}) \leq P_1(b \leftrightarrow a, b \leftrightarrow c_1, b \leftrightarrow c_2)\label{eq:locprobbound}.
    \end{equation}
    
    \begin{figure}[h!]
        \centering
        \includegraphics[scale=0.75]{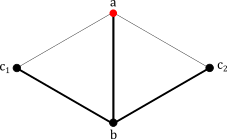}
    \end{figure}
    
    Depending on the type of vertices $a,b,c_1,c_2$ are, the pairs $(a,b),(b,c_1),(b,c_2)$ may be geometric or not.

    Suppose $m$ of these pairs are geometric (when $m=0$ none of them). Then at least $m$ of the vertices $b,c_1,c_2$ are type-B. These $m$ vertices can be chosen in at most $k$ different ways, and the remaining vertices in at most $n$ different ways. From Lemma \ref{lemma:marginalprob} we know that every geometric pair $i,j$ has connection probability $p_{ij} = \Theta(1 \wedge w_i w_j/k = O(w_i w_j/k)$. Then, from \eqref{eq:locprobbound} and applying Lemma \ref{lemma:geopathbound}, we have that the contribution to \eqref{eq:locvarbound} given by the combination of vertices such that $m$ among $b,c_1,c_2$ are type-B is
    \begin{equation}
    \begin{split}
        \frac{n^2}{w_a^4}\sum_{\substack{b,c_1,c_2 \\ m \text{ are type-B}}}\frac{P_1(\triangle_{a,b,c_1},\triangle_{a,b,c_2})}{w_{b}^2w_{c_1}w_{c_2}} &\leq \frac{n^2}{w_a^4}\sum_{\substack{b,c_1,c_2 \\ m \text{ are type-B}}}\frac{O(w_a w_b^3 w_{c_1} w_{c_2})}{k^m n^{3-m} w_{b}^2w_{c_1}w_{c_2}} \\
        &\leq \frac{n^2}{w_a^4}\binom{3}{m}O(w_a) = O(n^2/w_a^3).
    \end{split}
    \end{equation}
    Therefore, summing up all possible choices for $m=0,1,2,3$, we have that 
    \begin{equation}
        \frac{n^2}{w_a^4}\sum_{b,c_1,c_2}\frac{P_1(\triangle_{a,b,c_1},\triangle_{a,b,c_2})}{w_{b}^2w_{c_1}w_{c_2}} = O(n^2/w_a^3).
    \end{equation}
 
    \item If $|\{b_1,c_1\} \cap \{b_2,c_2\}| = 2$, without loss of generality we may assume $b_1 = b_2 \equiv b$ and $c_1 = c_2 \equiv c$.
    The numerator in the summation term of Equation \eqref{eq:locvarbound} can be further upper bounded by
    \begin{equation}
        P_1(\triangle_{a,b,c}) \leq P_1(a \leftrightarrow b, b \leftrightarrow c) \label{eq:locprobbound2}
    \end{equation}

    \begin{figure}[h!]
        \centering
        \includegraphics[scale=0.75]{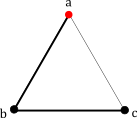}
    \end{figure}
    
    Similarly as before, depending on the type of vertices, the connection probabilities may be geometric or non-geometric.
    Let $m$ be the number of pairs among $(a,b),(b,c)$ which are geometric.
    As in the previous step, from \eqref{eq:locprobbound2} and applying Lemmas \ref{lemma:sumweights}-\ref{lemma:geopathbound} we have
    \begin{equation}
    \begin{split}
        \frac{n^2}{w_a^4}\sum_{b,c} \frac{P_1(\triangle_{a,b,c})}{w_{b}^2w_{c}^2} &\leq \frac{n^2}{w_a^4}\sum_{\substack{b,c \\ m \text{ are type-B}}}\frac{O(w_a w_b^2 w_{c})}{k^m n^{2-m} w_{b}^2w_{c}^2} \\
        &\leq \frac{n^2}{w_a^4}\binom{2}{m}O(w_a) = O(n^2/w_a^3).
    \end{split}
    \end{equation}
    Therefore, summing up all possible choices for $m=0,1,2$, we have that 
    \begin{equation}
        \frac{n^2}{w_a^4}\sum_{b,c}\frac{P_1(\triangle_{a,b,c})}{w_{b}^2w_{c}^2} = O(n^2/w_a^3).
    \end{equation}
\end{itemize}

Therefore, summing up all possible intersections for the sets $\{b_1,c_1\}, \{b_2,c_2\}$ we conclude that the the right hand side of Equation \eqref{eq:locvarbound} is $O(n^2/w_a^3)$.
\end{proof}

\subsection{Node-based detection test}

\begin{proof}[Proof of Theorem~\ref{thm:identify}]
Assume that $w_a \gg \log(n)$. 

We first calculate the type I error. When $a$ is type-A, from Proposition \ref{prop:expWiIRG} and Proposition \ref{prop:varianceWi} we have that $\Expa{\Wi(a)} = O(1), \Var(W(a)) = O(n^2/w_a^3)$. By Chebyshev's inequality,
\begin{align}\label{eq:nodetypeI}
    P\left(\Wi(a)>\frac{n}{w_a \sqrt{\log(n)}}\right) &\leq \frac{\Var(\Wi(a))}{\left( \frac{n}{w_a \sqrt{\log(n)}}(1+o(1)) \right)^2} = O\left(\frac{\log(n)}{w_a}\right).
\end{align}
Let $\tilde{n}$ denote the number of vertices with weight at least $t_n$, and $\tilde{k}$ the number of type B vertices with weight at least $t_n$. 
For the type II error, from Proposition \ref{prop:expWiGIRG} and Proposition \ref{prop:varianceWi} we have that $\Expa{\Wi(a)} = \Omega(n/w_a), \Var(W(a)) = O(n^2/w_a^3)$. for any type-B vertex $a$. Then, by Chebyshev's inequality
\begin{align}\label{eq:nodetypeII}
    P\left(\Wi(w_i)<\frac{n}{w_a \sqrt{\log(n)}}\right) &\leq \frac{\Var(W(a))}{\left(\frac{n}{w_a \sqrt{\log(n)}}-\Expa{\Wi(a)}\right)^2} \nonumber\\
    &= O\left( \frac{n^2/w_a^3}{(n/w_a)^2} \right) = O\left( \frac{1}{w_a} \right).
\end{align}
Thus, by~\eqref{eq:invdetailed} the expected number of misclassified type B vertices equals
\begin{equation}
    \tilde{k}K_2(h\log(n))^{-\tau}=o(\tilde{k})
\end{equation}
for some $K_2>0$. 


With weight threshold $t_n$ for which we apply the test, as $\tilde{k}$ is binomial with mean $nt_n^{1-\tau}$, Theorem A.1.4 in \cite{AlonSpencer} yields that for all $\varepsilon>0$
\begin{equation}
    \Prob{\mathcal{E}_1}:=\Prob{\tilde{k}<(1-\varepsilon)kt_n^{1-\tau}}\leq \exp(-kt_n^{1-\tau}((1-\varepsilon)\log(1-\varepsilon)+\varepsilon))= \exp(-\tilde{\varepsilon}kt_n^{1-\tau}),
\end{equation}
for some $\tilde{\varepsilon}>0$, and 
\begin{equation}
    \Prob{\mathcal{E}_2}:=\Prob{\tilde{n}<(1+\varepsilon)nt_n^{1-\tau}}\leq \exp(-nt_n^{1-\tau}((1+\varepsilon)\log(1-\varepsilon)-\varepsilon))= \exp(-\hat{\varepsilon}nt_n^{1-\tau}),
\end{equation}
for some $\hat{\varepsilon}>0$. 
Thus, as $k<n$,
\begin{equation}
    \Prob{\bar{\mathcal{E}_1}\cap \bar{\mathcal{E}_2}}\geq 1-\exp(-\zeta kt_n^{1-\tau})
\end{equation}
for some $\zeta>0$. 
By~\eqref{eq:nodetypeI} and~\eqref{eq:wumin1} with $w_0=t_n$, on $\bar{\mathcal{E}}_2$, the expected number of misclassified type A vertices equals
\begin{align}\label{eq:nodetypeIrisk}
     O\left(\sum_{a\in [n-k]:w_a>t_n}\frac{\log(n)}{w_a}\right) = O(\log(n)t_n^{-\tau}nt_n^{1-\tau}),
\end{align}
where $nt_n^{1-\tau}$ appears as on $\bar{\mathcal{E}}_2$, there are at most  $(1+\varepsilon) nt_n^{1-\tau}$ vertices of type A with weight at least $t_n$. Furthermore, misclassified type $B$ vertices are $o(\tilde{k})$ by~\eqref{eq:nodetypeII}.

Now,
\begin{align}
    \Exp{\frac{|\hat{V}_C^{t_n}\triangle V_C^{t_n}|}{2|V_C^{t_n}|}\mid V_C^{t_n}\geq 1}& = \Prob{\mathcal{E}_1 \cup \mathcal{E}_2} \Exp{\frac{|\hat{V}_C^{t_n}\triangle V_C^{t_n}|}{2|V_C^{t_n}|}\mid \mathcal{E}_1 \cup \mathcal{E}_2, V_C^{t_n}\geq 1} \nonumber\\
    & \quad +  \Prob{\bar{\mathcal{E}_1} \cap \bar{\mathcal{E}_2}} \Exp{\frac{|\hat{V}_C^{t_n}\triangle V_C^{t_n}|}{2|V_C^{t_n}|}\mid \bar{\mathcal{E}_1} \cap \bar{\mathcal{E}_2}}\nonumber\\
    & \leq n\exp(-\tilde{\varepsilon}kt_n^{1-\tau}) + \frac{ O(\log(n)t_n^{1-2\tau}n)}{(1-\varepsilon)kt_n^{1-\tau}}\nonumber\\
    & = o(1),
\end{align}
as long as 
\begin{equation}\label{eq:tncond}
     (n\log(n)/k)^{1/\tau} \ll t_n \ll k^{1/(\tau-1)}.
\end{equation}
\end{proof}

\appendix

\section{Expected degrees in $H_1$}\label{app:consistency_correction}
We show now why, under the alternative hypothesis $H_1$, a small correction on the connection probabilities has to be made. We first compute the connection probability of a vertex with weight $w_i$, then we apply it to find the expected degree of a vertex with weight $w_i$. 

In this Appendix, we call $\mathbb{E}_X$ the expected value over the random variable $X$.

\begin{lemma}\label{lemma:marginalprob}
Under $H_1$, any two type-B vertices $i$ and $j$ are connected with probability
\begin{equation}
    P_1(i \leftrightarrow j) = \mathbb{E}_{x_i,x_j}[p_{ij}(w_i,w_j,x_i,x_j)] = \begin{cases}
    \frac{1}{1 + C_1} & \text{ if $\frac{w_i w_j}{k} \geq \frac{\mu}{2^d}$}\\
    \frac{C_1}{1 + C_1} \frac{w_i w_j}{k} + o(\frac{w_i w_j}{k}) & \text{ otherwise.}\\
    \end{cases}
\end{equation}
\end{lemma}

\begin{proof}
First, observe that 
\begin{align}
    \mathbb{E}_{x_j}[p_{ij}(w_i,w_j,x_i,x_j)] &=  \frac{1}{1 + C_1}\int_{\mathbb{T}^d} \left(\frac{w_i w_j}{n\mu ||x_i-x_j||^d} \wedge 1 \right)^{\gamma} \dd x_j \nonumber\\
    &= \frac{1}{1 + C_1}\int_0^{1/2} \left(\frac{w_i w_j}{n\mu r^d} \wedge 1 \right)^{\gamma} \Prob{||x_j||=r} \dd r \nonumber\\
    &= \frac{1}{1 + C_1}\int_0^{1/2} \left(\frac{w_i w_j}{n\mu 2^{-d}V(r)} \wedge 1 \right)^{\gamma} V'(r) \dd r 
\end{align}
because $\bx_i$ is fixed, and $\bx_j$ is uniformly distributed, hence $V(r) = \Prob{||x_j|| \leq r} = (2r)^d$. Then, averaging over $x_i$, and applying integral substitution
\begin{equation}
\begin{split}
    \mathbb{E}_{\bx_i}\left[\mathbb{E}_{\bx_j}[p_{ij}(w_i,w_j,x_i,x_j)]\right] &= \frac{1}{1 + C_1} \int_{[0,1]^d} \left( \int_0^{1} \left(\frac{w_i w_j}{n\mu 2^{-d}V} \wedge 1 \right)^{\gamma} dV \right) \dd x_i \\
    &= \frac{1}{1 + C_1} \int_0^{1} \left(\frac{w_i w_j}{n\mu 2^{-d}V} \wedge 1 \right)^{\gamma} \dd V.
\end{split}
\end{equation}
We denote $r_0 = \frac{w_i w_j}{n \mu 2^{-d}}$. 
\begin{itemize}
    \item When $r_0 \geq 1$, then the latter integrand is always $1$, and $\mathbb{E}_{x_i,x_j}[p_{ij}(w_i,w_j,x_i,x_j)] = 1$.
    \item When $r_0 < 1$, 
    \begin{align}
        \mathbb{E}_{x_i,x_j}[p_{ij}(w_i,w_j,x_i,x_j)] &= \frac{1}{1 + C_1} \left(\int_0^{r_0} 1 dV + \int_{r_0}^{1} \left(\frac{r_0}{V} \right)^{\gamma} dV \right)\nonumber\\
        &=  \frac{1}{1 + C_1} \left(r_0 + \left[ \frac{r_0^{\gamma}}{\gamma -1}\left( r_0^{-\gamma + 1} - 1 \right) \right]\right) \nonumber\\
        &=  \frac{1}{1 + C_1}r_0 (1 + (\gamma -1)^{-1}) + O(r_0^{\gamma}). \label{marginal}
    \end{align}
\end{itemize} 
\end{proof}

\begin{lemma}\label{lemma:consistency_correction}
Under $H_1$, any vertex $i$ has expected degree $\mathbb{E}[d_i]=w_i(1+o(1)).$
\end{lemma}

\begin{proof}
The expected degree of any vertex $i$ is:
\begin{equation*}
    E_1[d_i] 
    = E_1\left[\sum_{j \in V} \mathbbm{1}_{\{i \sim j\}}\right]
\end{equation*}
\begin{itemize}
    \item If the vertex $i$ is type-A, it connects to any other vertex using \eqref{eq:IRG_connection_rule}.
    In this case, the expected degree of $i$ is upper and lower bounded by
    \begin{equation*}
        E_1\left[\sum_{j \in V} \mathbbm{1}_{\{i \sim j\}}\right] = \sum_{j \in V} p_{ij}(w_i,w_j) \leq \frac{w_i}{\mu n} \sum_{j \in V} w_j = w_i (1 + o(1))
    \end{equation*}
    and
    \begin{equation*}
        E_1\left[\sum_{j \in V} \mathbbm{1}_{\{i \sim j\}}\right] = \sum_{j \in V} p_{ij}(w_i,w_j) \geq \frac{w_i}{k} \sum_{j \in V \setminus [n-k]} w_j \mathbbm{1}_{\{w_j \leq \mu n / w_i\}} = w_i (1 + o(1))
    \end{equation*}
    Then, $E_1[d_i] = w_i (1 +o(1))$.
    
    \item If the vertex $i$ is type-B, it connects to a type-A vertex $j$ using \eqref{eq:corr_IRG_connection_rule} and to a type-B vertex $k$ using \eqref{eq:corr_GIRG_connection_rule}. In this case,
    \begin{equation*}
        E_1\left[\sum_{j \in V} \mathbbm{1}_{\{i \sim j\}}\right] = \sum_{j \in [n-k]} p_{ij}(w_i,w_j) + \sum_{j \in V \setminus [n-k]} \mathbb{E}_{x_i,x_j}[p_{ij}(w_i,w_j,x_i,x_k)] =: A + B.
    \end{equation*}
    From the same computations as in the previous case, we have that $A = \frac{n-k}{n} \frac{w_i}{1 + C_1}(1 + o(1)) = \frac{w_i}{1 + C_1}(1 + o(1))$. On the other hand, when $i,j$ are type-B, 
    \begin{equation*}
        P(i \leftrightarrow j) = \begin{cases}
        \frac{1}{1 + C_1} & \text{ if $\frac{w_i w_j}{k} \geq \frac{\mu}{2^d}$}\\
        \frac{C_1}{1 + C_1} \frac{w_i w_j}{k} + o(\frac{w_i w_j}{k}) & \text{ otherwise}\\
        \end{cases}
    \end{equation*}
    from Lemma \ref{lemma:marginalprob}. Thus, we have the upper and lower bounds
    \begin{align*}
        B &\leq \frac{C_1}{1 + C_1} \frac{w_i}{k} \sum_{j \in V \setminus [n-k]} w_j (1+o(1)) = \frac{C_1}{1 + C_1} w_i (1 + o(1)), \\
        B &\geq \frac{C_1}{1 + C_1} \frac{w_i}{k} \sum_{j \in V \setminus [n-k]} w_j (1+o(1)) \mathbbm{1}_{\{w_j \leq \mu k / w_i 2^d\}} = \frac{C_1}{1 + C_1} w_i (1 + o(1)).
    \end{align*}
    Then, $E_1[d_i] = A + B = w_i (1 + o(1))$.
\end{itemize}
\end{proof}

\section{Auxiliary results}

\begin{lemma}\label{lemma:geopathbound}
Given the weight sequence $(w_i)_{i \in V}$, suppose $G$ is a graph sampled under the alternative hypothesis $H_1$. Let $p \geq 2$ and $\mathcal{P} = (a_1,...,a_p)$ an ordered set of vertices in the graph $G$. Denote by $m = |\{i  : a_i,a_{i+1} \text{ are type-B}\}|$ the number of pairs in $\mathcal{P}$ that are connected according to the GIRG connection probability. Then, 
\begin{align}
    P_1(a_1 \leftrightarrow a_2 \leftrightarrow ... \leftrightarrow a_p) &= \prod_{i=1}^{p-1} P_1(a_i \leftrightarrow a_{i+1}) \label{eq:geopathind}\\ & = O\left(\frac{w_{a_1}^{-1} w_{a_p}^{-1} \prod_{i = 1}^p w_{a_i}^2}{k^m n^{p-m-1}}\right) \label{eq:geopathbound}
\end{align}
\end{lemma}
\begin{proof}
We first prove \eqref{eq:geopathind}, using an inductive argument on the size of the path $\mathcal{P}$. When $p=2$ the equality is trivial. Next, suppose the equality holds for $p=s$, and let us prove it for $p=s+1$. Consider the pair $(a_{s}, a_{s+1})$.\\
If at least one between $a_s$ and $a_{s+1}$ is a type-A vertex, then the random variable $\mathbbm{1}_{a_s \leftrightarrow a_{s+1}}$ is a Bernoulli random variable with parameter $\min(1, w_{a_s}w_{a_{s+1}}/\mu n)$, which is independent from the rest of the path. Therefore, $P_1(a_1 \leftrightarrow a_2 \leftrightarrow ... \leftrightarrow a_{s+1}) = \prod_{i=1}^{s-1} P_1(a_i \leftrightarrow a_{i+1}) \cdot P_1(a_s \leftrightarrow a_{s+1})$.\\
If both $a_s$ and $a_{s+1}$ are type-B, then the probability for the edge $(a_s, a_{s+1})$ depends on the positions $\bx_{a_s}$ and $\bx_{a_{s+1}}$. However, the probability that the entire path $(a_1,...,a_s)$ is present in the graph may depend on the positions $\bx_{a_1}...\bx_{a_s}$. Then, we may condition on the position $\bx_{a_s}$ and obtain
\begin{equation}
    P_1(a_1 \leftrightarrow a_2 \leftrightarrow ... \leftrightarrow a_{s+1}) = \int_{\mathbb{T}^d} P_1(a_1 \leftrightarrow a_2 \leftrightarrow ... \leftrightarrow a_{s} \;|\; \bx_{a_s}) P_1(a_s \leftrightarrow a_{s+1} \;|\; \bx_{a_s}) \dd P_1(\bx_{a_s} = x). 
\end{equation}
Now, by the symmetry of the infinity norm on the torus, follows that $P_1(a_s \leftrightarrow a_{s+1} \;|\; \bx_{a_s}) = P_1(a_s \leftrightarrow a_{s+1})$. Therefore,
\begin{align}
    P_1(a_1 \leftrightarrow a_2 \leftrightarrow ... \leftrightarrow a_{s+1}) &= P_1(a_s \leftrightarrow a_{s+1}) \int_{\mathbb{T}^d} P_1(a_1 \leftrightarrow a_2 \leftrightarrow ... \leftrightarrow a_{s} \;|\; \bx_{a_s}) \dd \mathbb{P}(\bx_{a_s} = x) \\
    &= P_1(a_s \leftrightarrow a_{s+1}) P_1(a_1 \leftrightarrow a_2 \leftrightarrow ... \leftrightarrow a_{s}), 
\end{align}
and the proof is concluded applying the inductive hypothesis. \\

Next, we prove \eqref{eq:geopathbound}. In the path $m$ pairs of vertices are connected according to the GIRG edge probability. Assume $(a_r,a_{r+1})$, for some $1 \leq r < p$ is one such a pair. Then, applying Lemma \ref{lemma:marginalprob} we have that $P_1(a_r \leftrightarrow a_{r+1}) = O(w_{a_{r}}w_{a_{r+1}}/k)$. The remaining $p-m-1$ pairs of vertices are connected following the IRG edge rule. In that case, if $(a_r,a_{r+1})$ is such a pair, for some $r$, we have the upper bound $P_1(a_r \leftrightarrow a_{r+1}) = \min(1, w_{a_{r}}w_{a_{r+1}}/\mu n) = O(w_{a_{r}}w_{a_{r+1}}/n)$. Combining these two facts together, we conclude that
\begin{equation}
    \prod_{i=1}^{p-1} P_1(a_i \leftrightarrow a_{i+1}) = O\left(\frac{w_{a_1}^{-1} w_{a_p}^{-1} \prod_{i = 1}^p w_{a_i}^2}{k^m n^{p-m-1}}\right)
\end{equation}
regardless of which of the $m$ pairs form between type-B vertices.
\end{proof}

\bibliographystyle{abbrv}

\bibliography{references}

\end{document}